\documentclass[12pt]{article}
\makeindex
\usepackage{epsfig,amsmath,a4wide,amssymb,amsfonts,amsbsy,xy,latexsym,epsf,pstricks,verbatim,times}
\xyoption{all}

\let\set\mathbb
\def\<<{\leavevmode
  \raise0.28ex\hbox{$\scriptscriptstyle\langle\!\langle$}\nobreak
  \hskip -.6pt plus.3pt minus.2pt\,}
\def\>>{\,\nobreak\hskip -.6pt plus.3pt minus.2pt
  \raise0.28ex\hbox{$\scriptscriptstyle\rangle\!\rangle$}}

\def\talpha{{\tilde \alpha}}
\def\tbeta{{\tilde \beta}}
\def\vp{{\varphi}}
\def\AA{{\set A}}

\def\hT{{\hat T}}

\def\Hom{\mathop{\rm{Hom}}\nolimits }
\def\IM{\mathop{\rm{Im}}\nolimits }

\def\CC{{\set C}}

\def\Div{{\rm Div}}

\def\End{\mathop{\rm End }}
\def\FF{{\set F}}

\def\FQ{{\FF _Q}}
\def\Fp{{\FF _p}}
\def\Fl{{\FF _\ell}}
\def\Fqb{{{\bar {\FF}} _q}}
\def\Fpb{{{\bar {\FF}} _p}}

\def\Fq{{\FF _q}}
\def\Fqs{{\FF ^*_q}}

\def\Gal{\mathop{\rm{Gal}}\nolimits }

\def\Ker{\mathop{\rm{Ker}}\nolimits }

\def\NN{{\set N}}
\def\PP{{\set P}}
\def\GG{{\set G}}

\def\II{{\set I}}
\def\PP{{\set P}}

\def\Pic{\mathop{\rm{Pic}}\nolimits }
\def\QQ{{\set Q}}
\def\Qb{{\bar \QQ}}

\def\TT{{\set T}}
\def\Id{{\rm{Id}}}

\def\GL{\mathop{\rm{GL}}\nolimits }

\def\ZZ{{\set Z}}
\def\agot{{\mathfrak a}}
\def\Bgot{{\mathfrak B}}
\def\Cgot{{\mathfrak C}}
\def\Dgot{{\mathfrak D}}

\def\bQ{{\bar \QQ}}
\def\bgot{{\mathfrak b}}
\def\cA{{\cal A}}
\def\cB{{\cal B}}
\def\cC{{\cal C}}
\def\cD{{\cal D}}

\def\cF{{\cal F}}

\def\cH{{\cal H}}

\def\cJ{{\cal J}}
\def\cK{{\cal K}}
\def\cL{{\cal L}}

\def\cM{{\cal M}}
\def\cN{{\cal N}}
\def\cO{{\cal O}}
\def\cP{{\cal P}}

\def\cR{{\cal R}}
\def\cS{{\cal S}}
\def\cT{{\cal T}}
\def\cU{{\cal U}}

\def\cX{{\cal X}}
\def\cY{{\cal Y}}

\def\mgot{{\mathfrak  m}}
\def\mmu{{\set \mu}}

\def\pgot{{\mathfrak p}}

\newtheorem{lemma}{Lemma}

\newtheorem{theorem}{Theorem}

\newtheorem{definition}{D{e}finition}
\newtheorem{remark}{Remark}

\begin{document}
\author{J.-M. Couveignes\thanks{Institut de  Math\'ematiques  de Toulouse,
Universit\'e de Toulouse et CNRS, D{\'e}partement de Math{\'e}matiques et
Informatique, Universit{\'e} Toulouse 2, 5 all{\'e}es Antonio Machado, 31058 Toulouse
c{\'e}dex 9.}}
\title{Linearizing  torsion  classes in the Picard group of
  algebraic curves over finite fields\thanks{Research supported by the Agence
    Nationale de la Recherche (projet blanc ALGOL).}}

\maketitle

\bibliographystyle{plain}

\begin{abstract}
We address the problem of computing 
in  the group of $\ell^k$-torsion 
rational points of the jacobian variety 
of algebraic curves  over  finite fields, with a view toward
 computing modular representations.
\end{abstract}
\tableofcontents

\section{Introduction}\label{section:introduction}

Let $\Fq$ be a finite field of characteristic $p$ 
and $\AA^2\subset \PP^2$ the affine and projective planes over $\Fq$
and
 $C\subset \PP^2$ a plane projective absolutely irreducible reduced
curve  over $\Fq$
and $\cX$ its smooth projective model and $\cJ$ the jacobian variety of
$\cX$. 
Let $g$ be the genus of $\cX$ and $d$ the degree of $C$.

We assume that we are given the numerator of the zeta function of the
function
field $\Fq(\cX)$. So we know the characteristic polynomial of the
Frobenius endomorphism $F_q$ of  $\cJ$. This is a monic degree $2g$ polynomial
$\chi(X)$ with integer coefficients.

Let $\ell\not = p $ be a prime integer and let $n=\ell^k$ be a power
of $\ell$. We
look for a {\it nice generating set} for the group $\cJ[\ell^k](\Fq)$ 
of  $\ell^k$-torsion points in $\cJ(\Fq)$. By {\it nice}  we mean that the
generating set $(g_i)_{1\le i\le I}$ should induce  a decomposition of
$\cJ[\ell^k](\Fq)$ as  a direct product $\prod_{1\le i\le I} <g_i>$
of cyclic subgroups with non-decreasing orders. 

Given such a generating set and an $\Fq$-endomorphism of $\cJ$,  we
also want to describe the action
of this endomorphism on  $\cJ[\ell^k](\Fq)$ by an $I\times I$ integer matrix.

In section \ref{section:classics} we recall how to compute in the 
Picard group
$\cJ(\Fq)$. Section \ref{section:picking}  gives a naive algorithm for picking random elements in this
group. Pairings are useful when looking for relations between
divisor classes. So we recall how to compute pairings in section
\ref{section:pairings}. Section \ref{section:divisiblegroups} is
concerned with characteristic subspaces for the action of Frobenius
inside the $\ell^\infty$-torsion of $\cJ(\Fqb)$. In section
\ref{section:kummer} we look for  a convenient surjection from $\cJ(\Fq)$ onto
its
$\ell^k$-torsion subgroup. We use the Kummer exact sequence and the
structure
of the ring generated by the Frobenius endomorphism. In section \ref{section:relations} we give
an  algorithm that, on input a degree $d$ plane
projective curve  over $\Fq$, plus some information on its
singularities, and the zeta function of its function field, 
   returns a nice generating set
for the group of $\ell^k$-torsion points inside $\cJ(\Fq)$ in
probabilistic  polynomial time in $\log q$, $d$ and $\ell^k$. 
Sections 
\ref{section:modularcurves} and \ref{section:modularcurves2} are devoted to two families of modular curves. We
give
a nice plane model for such curves. The general algorithms presented
in
section \ref{section:relations} are then applied to these modular curves in
section
\ref{section:ramanujan} in order to compute explicitly the modular
representation modulo $\ell$ associated with  the discriminant modular form (level $1$
and weight $12$).  
This modulo $\ell$ representation $V_\ell$ is seen  as a
subgroup of order $\ell^2$ inside the $\ell$-torsion of $J_1(\ell)/\QQ$. The idea is to
compute the reduction modulo $p$ of the
group scheme  $V_\ell$ as a subgroup of   $J_1(\ell)/\Fp$,    for many small primes
$p$. One then lifts using the Chinese Remainder Theorem. 
This makes a connection with  Edixhoven's
program for  computing coefficients of modular forms.  My contribution to this
program is sketched in section \ref{section:prog}.  See \cite{edix2, arxiv}.
The core of Edixhoven's program is that if one knows $V_\ell$,  one can
efficiently compute the Ramanujan function $\tau(P)$ modulo
$\ell$ for  a large prime $P$. If we have enough primes $\ell$, we can
deduce the actual value of $\tau(P)$.

The last three sections present variants of
the main algorithm and auxiliary results.
Section \ref{section:semi} presents a simpler variant of the method of
section \ref{section:ramanujan}, that is particularly useful when the action of
the $p$-Frobenius 
on $V_\ell$  modulo $p$ is semisimple non-scalar. In the non-semisimple case,
this simpler method may  only produce a non-trivial subspace inside $V_\ell$
  modulo $p$.
 Section \ref{section:elem}
proves that this semisimplicity  condition holds  quite often indeed,
as expected. As a consequence,  one may compute the  representation $V_\ell$
associated with the discriminant form for at least half (say) the
primes $\ell$, using this simplified algorithm. This suffices for the purpose of computing the
Ramanujan function $\tau(P)$ at a large prime $P$ since we may afford
to skip half the auxiliary primes $\ell$. On the other hand, if one
wishes to compute a   representation modulo $\ell$ for a given
$\ell$,  then one should
be ready to face (at least theoretically) the  case when no small
prime $p$ is  semisimple for $\ell$. In that situation,
the 
simplified algorithm would only give 
a non-trivial subspace of $V_\ell$ modulo $p$ for many primes $p$.

Section
\ref{section:lift} addresses  the problem of computing $V_\ell$ from
all the knowledge we have collected concerning $V_\ell \bmod p$ for many
small primes $p$.
It requires  a sort of interpolation  theorem in the context
of polynomials with integer coefficients. The goal is to recover a
polynomial $P(X)$ once given a collection of non-trivial factors of
$P(X) \bmod p$ for many primes $p$.  This helps recovering 
$V_\ell/\QQ$  once given a subspace in  its reduction modulo
$p$ for enough small primes $p$.

Altogether, this proves that the simplified algorithm, despite
the possibility of 
many  non-semisimple primes $p$, suffices to compute
$V_\ell/\QQ$ for all $\ell$.

\begin{remark} The symbol $\cO$ in this article stands for
a positive effective absolute constant. So any statement containing
this symbol becomes true if the symbol is replaced in every occurrence by some 
large enough real number.
\end{remark}

\begin{remark} By an algorithm in this paper we usually mean a
probabilistic (Las Vegas) algorithm. This is an algorithm that
succeeds with probability $\ge \frac{1}{2}$. When it fails, it gives no
answer.
 In some places we shall give
deterministic algorithms or probabilistic (Monte-Carlo) algorithms, but
this will be stated explicitly. A Monte-Carlo algorithm gives a
correct answer with probability $\ge \frac{1}{2}$. But it may give an
incorrect answer with probability $\le \frac{1}{2}$. A Monte-Carlo
algorithm can be turned into a Las Vegas one, provided we can
efficiently check the correctness of the result.  One reason for using probabilistic
Turing machines is that in many places it will be necessary (or at
least wiser) to decompose a divisor as a sum of places. This is the
case in particular for the conductor of some plane curve. Another
more intrinsically probabilistic algorithm  in this paper is the one that searches
for generators of the Picard group. 
\end{remark}

\section{Context: the inverse Jacobi problem}\label{section:prog}

The  initial motivation for this work 
is a discussion I had in 2000 with Bas Edixhoven about 
his program 
aiming at polynomial time computation of 
coefficients of modular forms.

He asked  how one can 
compute  (e.g.) the decomposition field of 
the dimension two modulo $\ell$ Galois representation  $V_\ell$ 
associated to the discriminant modular form $\Delta$.
This amounts to computing the field of moduli of some 
very special $\ell$-cyclic  
coverings of $X_1(\ell)$.

I had some  experience in explicit computation of coverings 
using numerical techniques and got the impression that 
a purely algebraic approach would fail to solve 
such a problem. This is because $V_\ell$, however small it is,
is lost in the middle of the full $\ell$-torsion
of $J_1(\ell)$. And the latter is a huge dimension zero
variety (its number of geometric points
 is exponential in $\ell$).

The second time I discussed this question with Edixhoven, it became clear
that we had two options. We might compute $V_\ell$ inside
the complex torus of $J_1(\ell)$ and evaluate a theta function
at some point $x$ in $V_\ell$. Edixhoven convinced me that this approach
was unlikely to succeed since the number of terms
to be considered in the expansion of the theta function
would be exponential in $\ell$, even for a poor
accuracy.
Another possibility  was to solve the inverse Jacobi problem
for $x$ and find a divisor $D=P_1+\dots+P_g-gO$  in the 
 class associated to $x$ in the Picard group of $X_1(\ell)$.
Then one would pick a function $f$ on $X_1(\ell)$ and
evaluate $F(x)=f(P_1)+\dots +f(P_g)$ for example.

Solving the  inverse Jacobi problem seemed easy.
Indeed one could pick any divisor $D^0=P_1^0+\dots+P_g^0-gO$
of the above form on $X_1(\ell)$ and compute its image
$x^0$ by the Jacobi map. Then one would move slowly from
$x^0$ to $x$ inside the complex torus $J_1(\ell)(\CC)$. At each step the corresponding divisor 
would be computed from the previous one using Newton's method.

Although the Jacobi map is birational, it is not quite an isomorphism
however. It has a singular locus and it was not clear  how
one could avoid this obstacle  in the journey from $x^0$ to $x$.

It was decided that I would think about how to solve
this problem while Edixhoven  would prove good bounds on the height of
the algebraic number $F(x)$ coming out of the algorithm.
Edixhoven   first proved the analogous bound in the function field case.
Then,  Bas  Edixhoven and  Robin de Jong, using Arakelov theory
and results by Merkl in \cite{arxiv} or J. Jorgenson and  J. Kramer in \cite{JJ}, proved the bound for
the height of $F(x)$.

On my side, I was trying to avoid the singular locus.
I believe that in general, the problem of avoiding the singular locus
might  very well be NP-complete. Indeed, if the curve under consideration
is very close to the boundary of the moduli space, the problem
takes  a discrete  aspect: the curve has long tubes and sometimes
one may have to decide to push one point through one tube
or the other one. In case one makes the wrong decision, one may
be lost for ever. The problem can be phrased in a more
mathematical way: if the curve is (close to) a Mumford curve,
solving the inverse Jacobi problem assumes one can solve the 
discrete counterpart for it: solving the Jacobi problem
for a finite graph; namely the intersection graph of the curve.
See \cite{boundary}  theorem Theorem 2.1 and the following remark
for a statement of this problem, that I suspect is very hard when
the genus of the graph tends to infinity.

Of course one may expect that $J_1(\ell)$ keeps far enough 
from the boundary of its moduli space when $\ell$
tends to infinity. However, I was not able to give a proof
that the above ideas do succeed in solving the inverse
Jacobi problem, even  for these curves.  I had to build on a rather
different idea and proved in \cite{jacojaco} that for
$X_0(\ell)$ at least, solving the inverse Jacobi problem
is deterministic polynomial time in $\ell$ and the required 
precision. 

The first version of \cite{jacojaco} was
ready in January 2004.
Extending this result to any modular curve is just
a technical problem, but I confess I was tired with technicalities
and I stopped there with  the complex method. 

Starting in August 2003 I decided to look for a $p$-adic 
analogue of this complex
method: looking for a $p$-adic  approximation instead of a 
complex one.
After some hesitation I realized that computing modulo several
small primes $p$ and then lifting using the Chinese remainder 
would lead to a simpler algorithm.  This text  gathers  the
results of this research.
The methods presented here are  the discrete counterpart of the
ones in \cite{jacojaco}. The essence of theorem
\ref{theorem:computingltorsionmodp}
is that the discrete method presented in this paper
applies to modular curves $X_1(\ell)$.
This is exactly what is needed for the purpose of computing the Ramanujan
function.

 The complex
approach is more tedious but leads to deterministic  algorithms. The
main reason is that the set of complex points in the jacobian is a
connected
topological space. The modulo $p$ approach that we present
here seems intrinsically
probabilistic, because one has to find generators of Picard
groups of curves over finite fields.

I should also say that the complex approach was not 
abandoned  since Johan Bosman started in June 2004 his PhD with
Edixhoven on this topic and he succeeded in explicitly  computing some  $V_\ell$
using the complex method. See \cite{bosman}. He built on  the Newton
approach to solving the inverse Jacobi problem, as sketched above.
This shows that the singular locus of the Jacobi map
is not so disturbing after all, at least in practice.

Several  sections in this text  have been included in Edixhoven's report
\cite{arxiv}. Many thanks are due to Bas Edixhoven and Robin de Jong 
for useful
discussions, suggestions,  and comments. 

Many thanks also 
to John Cremona and the
anonymous referee for reading in detail this long 
manuscript and for their  useful comments.

\section{Basic algorithms for plane curves}\label{section:classics}

We  recall elementary results about computing 
in the Picard group of an algebraic curve over a finite field. See
\cite{hache, volcheck}.

\subsection{Finite fields}
We should first explain how finite fields are represented.
The base field $\Fq$ is given by an irreducible polynomial $f(X)$  with degree
$a$
and coefficients in $\Fp$ where $p$ is the characteristic and
$q=p^a$. So $\Fq$ is $\Fp[X]/f(X)$. An extension of $\Fq$ is given
similarly by an irreducible polynomial in $\Fq[X]$. Polynomial
factoring in $\Fq[X]$
is probabilistic polynomial time in $\log q$ and the degree of the polynomial to be
factored.  

\subsection{Plane projective curves and their smooth model}

We now explain how curves are supposed to be represented in
this
paper. 

To start with, a  projective  plane curve $C$ over $\Fq$ is  given by  a degree $d$ homogeneous
polynomial $E(X,Y,Z)$ in the three
variables $X$, $Y$ and $Z$, with coefficients in $\Fq$.
The curve $C$ is assumed to be absolutely irreducible and reduced.
By a {\it point} on $C$ we mean a geometric point (an element
of $C(\Fqb)$).
Any  $\Fqb$-point on $C$ can be represented by its affine or projective
coordinates. 

Let $\cX$ be a smooth model of $C$.
There is a desingularization map  $\cX\rightarrow C$.
If 
$P\in \cX(\Fqb)$ is  a geometric point on $\cX$
 above a singular point $S$
on $C$, we say that $P$ is a {\it singular branch}.

The {\it conductor} $\Cgot$ is an effective  divisor
on $\cX$ with even coefficients. 
Some authors call it the adjunction divisor.
Its support is made
of all singular branches. The conductor  expresses
the local behaviour of the map
$\cX\rightarrow C$.
See \cite[IV.1]{Serre}, \cite{gorenstein}.
We have  $\deg(\Cgot)=2\delta$ where
$\delta$
is the difference between the arithmetic genus $\frac{(d-1)(d-2)}{2}$ 
of $C$ and
the  geometric genus $g$ of $\cX$. Since $\delta\le \frac{(d-1)(d-2)}{2}$, 
the support of $\Cgot$
contains  at most  $\frac{(d-1)(d-2)}{2}$ geometric points 
in $\cX(\Fqb)$. So the field of definition of any
singular branch on $\cX$  is an extension of $\Fq$ with 
degree $\le  \frac{(d-1)(d-2)}{2}$.
A modern reference for singularities of plane curves
is \cite{casas} and especially section 5.8.

The smooth model $\cX$ of $C$ is not given as a projective variety.
Indeed, we shall only need  a nice local
description of $\cX$ above
 every singularity of $C$. This means we need a 
list of all singular points on $C$, and a 
list (a labelling) of all  points
in $\cX(\Fqb)$ lying  above every singularity of $C$ (the 
singular branches), and  a uniformizing
parameter at every such branch. We also need the  Laurent series expansions of
affine  plane coordinates in terms of all these uniformizing
parameters.


More precisely,  let
$P\in \cX(\Fqb)$ be a geometric point above a singular point $S$, and 
let $v$ be the corresponding
valuation.  The field of definition of $P$
 is an extension
field $\FF_P$ of $\Fq$ with  degree  $\le \frac{(d-1)(d-2)}{2}$. 
Let $x$ and $y$
be affine coordinates that vanish at the  singular point $S$ on $C$. We need
a local parameter $t$ at $P$ and  expansions $x=\sum_{k\ge
  v(x)}a_kt^k$ and $y=\sum_{k\ge
  v(y)}b_kt^k$ with coefficients in $\FF_P$. 

Because these expansions are not finite, we just assume we are given
an oracle that on input a positive integer $n$
returns  the first $n$ terms in all these expansions.


This is what we mean when we say the smooth model $\cX$ is given.

We may also assume that  we are given the conductor $\Cgot$ of $C$ as a
combination of singular branches with even  coefficients.  
The following algorithms still work if  the conductor
is  replaced by any  divisor
$\Dgot$  that 
is greater than the conductor and has polynomial degree  in $d$. Such
a divisor can be found easily: the singular
branches on $\cX$ are supposed to be known already, and the multiplicities
are bounded above by $\frac{(d-1)(d-2)}{2}$.

There are many families of curves for which such a smooth model
 can be
given as  a Turing machine that answers in probabilistic polynomial time in the
size $\log q$ of the field and  the degree $d$ of $C$ and the number
$n$ of requested significant terms in the parametrizations of singular branches. This is the case for
curves with ordinary multiple points for example. We shall show in
sections \ref{section:modularcurves} and \ref{section:modularcurves2}
that this is also the case for two  nice families of   modular curves.

\subsection{Divisors, forms,  and functions}

Smooth $\Fqb$-points on $C$ are  represented by their affine or projective
coordinates. Labelling for the  branches above singular points is
given in the description of $\cX$. 
So we know how to represent divisors on $\cX$.

For any   integer $h\ge 0$   
we set 

$$\cS_h=H^0(\PP^2/\Fq,\cO_{\PP^2/\Fq}(h))$$
\noindent  the
$\Fq$-linear space of degree $h$ homogeneous polynomials in $X$,
$Y$, and $Z$.
It is a vector space of  dimension  $\frac{(h+1)(h+2)}{2}$ over $\Fq$. A
basis for it is made of all monomials of the form $X^aY^bZ^c$
with $a,b,c \in \NN$ and  $a+b+c=h$. 

We denote by

$$\cH_h=H^0(\cX/\Fq,\cO_{\cX/\Fq}(h))$$
\noindent  the space of forms of degree $h$ on $\cX$. Here 
$\cO_{\cX/\Fq}(h)$ is the pullback of $\cO_{\PP^2/\Fq}(h)$ to $\cX$. 

Let  $W$ be a degree $h$   form on $\PP^2$ having non-zero pullback  $W_\cX$ on  $\cX$. Let $H=(W_\cX)$ be the divisor of this restriction. The map $f\mapsto
\frac{f}{W_\cX}$ is a bijection from $H^0(\cX/\Fq,\cO_{\cX/\Fq}(h))$ to the linear space $\cL(H)$.

If $\Delta$ is  a divisor on $\cX$  we note $\cH_h(-\Delta)$ the
subspace of forms
in $\cH_h$  with
divisor $\ge \Delta$.
The dimension of $\cH_h(-\Cgot)$ is at least
$dh+1-g-\deg({\Cgot})$ and is equal to this number when it
exceeds $g-1$. This is the case if  $h\ge d$. The
dimension of 
$\cH_h(-\Cgot)$ is     greater than $2g$ if $h\ge 2d$.

 The image of the restriction map $\rho : \cS_h\rightarrow \cH_h$
contains   $\cH_h(-\Cgot)$  according to Noether's residue theorem \cite[Theorem 7]{gorenstein}.

We set   $S_C=\cS_{2d}$ and  $\cH_C=\cH_{2d}(-\Cgot)$,
and $H_C=\rho^{-1}(\cH_C)\subset S_C$ and
$K_C=\Ker(\rho)\subset H_C$.
 So we have $0\rightarrow K_C
\rightarrow H_C \rightarrow \cH_C \rightarrow 0$.

To find linear equations for  $H_C \subset S_C$ we consider a
generic homogeneous form
$F(X,Y,Z)=\sum_{a+b+c=2d}\epsilon_{a,b,c}X^aY^bZ^c$ of degree $2d$
in  $X$, $Y$  and $Z$.
For every branch $P$  above a singular point $S\in C$ 
(assuming for example  that  $S$ has non-zero
$Z$-coordinate)
we replace
in $F(\frac{X}{Z},\frac{Y}{Z},1)$  the affine coordinates 
$x=\frac{X}{Z}$ and $y=\frac{Y}{Z}$
by their   expansions as  series in  the local parameter $t_P$ at this
branch. We ask  the resulting series in $t_P$ to have
valuation at least the multiplicity of $P$ 
in the conductor $\Cgot$.
Every singular branch thus  produces linear equations
in the $\epsilon_{a,b,c}$.  The
collection of all such equations defines the subspace $H_C$.

A basis for the  subspace $K_C\subset H_C\subset S_C$
consists of all  $X^aY^bZ^cE(X,Y,Z)$ with 
$a+b+c=d$. We fix a supplementary space $M_C$ to $K_C$
in $H_C$ and assimilate $\cH_C$ to it.

Given
a homogeneous form in three variables  one can  compute its  divisor on $\cX$ 
using resultants and the given expansions of affine coordinates in terms of
the
local parameters at every singular branch.
A function is given as a quotient of two forms.

\subsection{The Brill-Noether algorithm}\label{subsubsection:brill}

Linear spaces of forms computed in the previous paragraph
allow us to  compute in the group $\cJ (\Fq)$ of $\Fq$-points
in the jacobian of $\cX$. We fix an effective $\Fq$-divisor
$\omega$ with
degree
$g$ on $\cX$. This $\omega$ will serve as an origin:
a point  $\alpha \in \cJ (\Fq)$ is
represented by a divisor $A-\omega$ in the corresponding linear
equivalence
class, where $A$ is an effective $\Fq$-divisor with degree $g$.
Given another point  $\beta \in \cJ (\Fq)$ by a similar divisor
$B-\omega$, we can compute the space $\cH_{2d}(-\Cgot-A-B)$ which is
non-trivial  and pick a non-zero form $f_1$ in it.  The divisor of $f_1$
is
$(f_1)=A+B+\Cgot+R$ where $R$ is an effective divisor with  degree
$2d^2-2g-2\delta$. The linear space $\cH_{2d}(-\Cgot-R-\omega)$
has dimension at least $1$. We pick a non-zero form $f_2$ in it. It
has divisor $(f_2)=\Cgot+R+\omega+D$ where $D$ is effective with
degree $g$. And $D-\omega$ is linearly equivalent to
$A-\omega+B-\omega$.

In order to invert the class $\alpha$  of $A-\omega$ we pick a
non-zero form  $f_1$ in $\cH_{2d}(-\Cgot-2\omega)$.  The divisor of $f_1$
is
$(f_1)=2\omega+\Cgot+R$ where $R$ is an effective divisor with  degree
$2d^2-2g-2\delta$. The linear space $\cH_{2d}(-\Cgot-R-A)$
has dimension at least $1$. We pick a non-zero form $f_2$ in it. It
has divisor $(f_2)=\Cgot+R+A+B$ where $B$ is effective with
degree $g$. And $B-\omega$ is linearly equivalent to
$-(A-\omega)$.

This algorithm works just as well if we replace $\Cgot$ by some
$\Dgot \ge \Cgot$ having polynomial degree  in $d$. 

\begin{lemma}[Arithmetic operations in the jacobian]\label{lemma:arithmeticoperations}
Let $C/\Fq$ be a degree $d$ plane projective absolutely irreducible reduced
curve. Let $g$ be the geometric genus of $C$.
Assume we are given the smooth model $\cX$ of $C$ and a 
   $\Fq$-divisor with degree
$g$  on $\cX$, denoted $\omega$. We assume $\omega$
is given as a difference between two effective divisors with degrees bounded
by a polynomial in $d$. This $\omega$ serves as an origin.
Arithmetic operations in the Picard group  $\Pic^0(\cX /\Fq)$ can be performed
in   time  polynomial  in $\log q$ and $d$. This includes
addition, substraction and comparison of divisor classes.
\end{lemma}
  
If $\omega$ is not effective, we use lemma \ref{lemma:explicitRR} below to
compute a non-zero function $f$  in $\cL(\omega)$ and we
write $\omega'=(f)+\omega$. This is an effective divisor with degree
$g$. We replace $\omega$ by $\omega'$ and finish as in the paragraph
before lemma \ref{lemma:arithmeticoperations}
\hfill $\Box$

We now recall
the principle  of the Brill-Noether algorithm for computing
complete linear series. Functions in $\Fq(\cX)$ are represented as
quotients of forms.

\begin{lemma}[Brill-Noether]\label{lemma:explicitRR}
There exists an algorithm that on input 
 a degree $d$ plane projective absolutely irreducible reduced
curve $C/\Fq$ and the  smooth model $\cX$ of $C$
 and  two effective $\Fq$-divisors $A$ and $B$ on
 $\cX$,  computes a basis for $\cL(A-B)$  in time  polynomial 
in $d$ and $\log q$ and the degrees of $A$ and $B$.
\end{lemma}

We assume $\deg(A)\ge \deg(B)$, otherwise $\cL(A-B)=0$. Let
$a$ be the degree of $A$. We let $h$ be the smallest integer
such that $h\ge 2d$ and
$hd+g+1>a+(d-1)(d-2)$. 

So the space $\cH_h(-\Cgot-A)$ is non-zero. It
 is contained in the 
image of the restriction map $\rho : \cS_h\rightarrow \cH_h$
 so that we can represent it
  as a subspace of $\cS_h$. We pick a non-zero  form
$f$  in   $\cH_h(-\Cgot-A)$ and compute its divisor
$(f)=\Cgot+A+D$.   

The space $\cH_h(-\Cgot-B-D)$  is contained in the 
image of the restriction map $\rho : \cS_h\rightarrow \cH_h$
 so that we can represent it
  as a subspace of $\cS_h$. We compute forms $\gamma_1$, $\gamma_2$,
  ..., $\gamma_k$ in $\cS_h$ such that their images by $\rho$ provide
  a basis for  $\cH_h(-\Cgot-B-D)$. A basis for $\cL(A-B)$ is made
  of the functions $\frac{\gamma_1}{f}$, $\frac{\gamma_2}{f}$, ...,
  $\frac{\gamma_k}{f}$. 
Again this algorithm works just as well if we replace $\Cgot$ by some
$\Dgot \ge \Cgot$ having polynomial degree  in $d$. 
\hfill $\Box$

We deduce an explicit moving lemma for divisors.

\begin{lemma}[Moving divisor lemma I]\label{lemma:moving}
There exists an algorithm that on input 
 a degree $d$ plane projective absolutely irreducible reduced
curve $C/\Fq$ and the  smooth model $\cX$ of $C$
 and  a degree zero  $\Fq$-divisor $D=D^+-D^-$ and an
effective divisor $A$ with degree $<q$ on 
 $\cX$ computes a divisor $E=E^+-E^-$ linearly equivalent to $D$ and
disjoint to $A$   in time  polynomial 
in $d$ and $\log q$ and the degrees of $D^+$,  and
$A$. Further
the degree of $E^+$ and $E^-$ can be taken to be
$\le 2gd$.
\end{lemma}

Let $O$ be an $\Fq$-rational divisor on $\cX$ such that
$1\le \deg (O)\le d$
 and disjoint to $A$.  We may take $O$ to be
a well chosen fiber of some plane coordinate function on $\cX$.
We  compute the linear space
$\cL=\cL(D^+-D^-+2g O)$. The subset of functions $f$ in $\cL$ such that
$(f)+D^+-D^-+2g O$ is not disjoint to $A$ is contained  in a union of
at most $\deg(A) < q$ hyperplanes. We conclude invoking lemma
\ref{lemma:inequ} below.\hfill $\Box$

There remains to state and prove the

\begin{lemma}[Solving inequalities]\label{lemma:inequ}
Let $q$ be  a prime power, $d\ge 2$  and $n\ge  1$ two  integers and
let $H_1$, ..., $H_n$ be hyperplanes inside  $V=\FF_q^d$,
each  given by a linear equation. Assume $n<q$. There exists a deterministic
algorithm
that finds a vector in $U=V-\bigcup_{1\le k\le n}H_k$ in time polynomial in $\log q$, $d$ and $n$.
\end{lemma}

This is proved by lowering  the dimension $d$. For $d=2$ we pick
any affine line $L$ in $V$ not containing the origin. We  observe that
there are at least $q-n$ points in 
$U\cap L=L-\bigcup_{1\le k\le n}L\cap H_k$. We enumerate points in $L$ until we
find one which is not in any $H_k$. This requires at most $n+1$
trials.

Assume now $d$ is bigger than $2$. Hyperplanes in $V$ are parametrized
by the projective space $\PP(\hat V)$ where $\hat V$ is the dual of
$V$. We enumerate points in $\PP(\hat V)$ until we find a hyperplane
$K$ distinct from every  $H_k$. We compute a basis for $K$ and an
equation for every $H_k\cap K$ in this basis. This way,  we have lowered
the dimension by $1$.\hfill $\Box$

We can strengthen a bit the moving divisor algorithm by removing the
condition that $A$ has degree $<q$.  Indeed, in case this condition is
not met, we call
$\alpha$ the smallest integer such that
$q^\alpha > \deg (A)$ and we set $\beta =\alpha+1$. We apply lemma
\ref{lemma:moving} after base change to the field with $q^\alpha$
elements and find a divisor $E_\alpha$. We call $e_\alpha$ the norm of
$E_\alpha$
from $\FF_{q^\alpha}$ to $\FF_q$. It is equivalent to $\alpha D$. We
similarly construct a divisor $e_\beta$ that is equivalent to $(\alpha+1)D$.
We
return the divisor $E=e_\beta-e_\alpha$. We observe that we can 
take $\alpha \le
1+\log_q \deg (A)$ so the degree of 
the positive part $E^+$ of $E$ is $\le 6g d(\log_q(\deg (A))+1)$.

\begin{lemma}[Moving divisor lemma II]\label{lemma:moving2}
There exists an algorithm that on input 
 a degree $d$ plane projective absolutely irreducible
curve $C/\Fq$ and the  smooth model $\cX$ of $C$
 and  a degree zero  $\Fq$-divisor $D=D^+-D^-$ and an
effective divisor $A$ on 
 $\cX$ computes a divisor $E=E^+-E^-$ linearly equivalent to $D$ and
disjoint to $A$   in time  polynomial 
in $d$ and $\log q$ and the degrees of $D^+$,  and
$A$. Further
the degree of $E^+$ and $E^-$ can be taken to be
$\le 6gd (\log_q(\deg (A))+1)$.
\end{lemma}

\section{A first approach to picking random divisors}\label{section:picking}

Given a finite field $\Fq$ and a plane projective absolutely
irreducible reduced curve $C$ over $\Fq$  with projective smooth  model
 $\cX$, we call $\cJ$ the jacobian of $\cX$
and we consider  
two related problems:  picking a random element 
in  $\cJ (\Fq)$  with (close to) uniform distribution and finding a
generating set for (a large subgroup of) $\cJ (\Fq)$. 
Let $g$ be the genus of $\cX$.
We  assume we are given  a degree $1$ divisor $O=O^+-O^-$ where $O^+$
and $O^-$ are effective, $\Fq$-rational and have degree bounded by an
absolute 
constant times $g$.

We know from \cite[Theorem 2]{MST} that   
 the   group $\Pic^0(\cX/\Fq)$ is generated by the classes 
$[\pgot -\deg(\pgot)O]$ where $\pgot$ runs over the set
of prime divisors of
degree $\le 1+2\log_q(4g-2)$.  For the convenience of the reader 
we quote this result as a lemma.

\begin{lemma}[M\"uller, Stein, Thiel]\label{lemma:Stein}
Let $K$ be an algebraic  function field of one variable
over $\Fq$. Let  $N\ge 0$ be an integer.
Let $g$ be the genus of $K$.
Let $\chi : \Div (K)\rightarrow \CC^*$ be a  character
of finite order which is non-trivial when restricted to
$\Div^0$. Assume  
that $\chi(\Bgot)=1$ for every prime divisor $\Bgot$ of degree
$\le N$. Then 

$$N<{2\log_q(4g-2)}.$$
\end{lemma}

If  $q<4g^2$,  the number of prime divisors of degree $\le 1+2\log_q(4g-2)$
is bounded by $\cO g^\cO$. So we  can compute easily
 a small generating set
for $\cJ(\Fq)$.

In the rest of this section, we will assume that the size 
 $q$ of the
field is greater than or equal to $4g^2$. This condition
ensures the
existence of a $\Fq$-rational point.

Picking efficiently and provably random elements in $\cJ (\Fq)$ with
uniform distribution  seems
difficult to us. We first give here an algorithm for efficiently constructing
random
divisors with a distribution that  is far from uniform but still
sufficient to construct a generating set for a large subgroup of
$\cJ (\Fq)$. Once given generators, picking random elements becomes
much easier.

Let $r$ be the smallest prime integer bigger than $30$, $2g -2$ and $d$.
We observe $r$ is less than $\max (4g-4,2d,60)$.

The set  $\cP(r,q)$ of $\Fq$-places with  degree $r$ 
on $\cX$ has cardinality

$$\#\cP(r,q)= \frac{\#\cX (\FF_{q^{r}}) -
  \#\cX (\FF_{q})}{r}.$$

So

$$(1-10^{-2})\frac{q^{r}}{r}\le \# \cP(r,q)  \le (1+10^{-2})\frac{q^{r}}{r}.$$

Indeed, $\left|\# \cX (\FF_{q^{r}})-q^{r}-1 \right|\le
2g q^{\frac{r}{2}}$ and $\left|\# \cX (\FF_{q})-q-1 \right|\le
2g q^{\frac{1}{2}}$.

 So   $\left | \# \cP(r,q) -
\frac{q^{r}}{r}\right| \le
\frac{4g+3}{r}q^{\frac{r}{2}}\le 8q^{\frac{r}{2}}$
and $8r q^{\frac{-r}{2}}\le r2^{3-\frac{r}{2}}\le
10^{-2}$ since  $r \ge 31$.

Since we are given a degree $d$  plane model $C$ for the curve
$\cX$,
we have a degree $d$ map $x : \cX \rightarrow \PP^1$.
Since $d<r$, the  function $x$ maps $\cP(r,q)$
to the set  $\cU(r,q)$ 
of monic prime polynomials of degree $r$ over
$\Fq$. The cardinality of   $\cU(r,q)$ is
$\frac{q^{r}-q}{r}$ so 

$$(1-10^{-9})\frac{q^{r}}{r}\le \# \cU(r ,q)  \le
\frac{q^{r}}{r}.$$

The fibers of the map $x : \cP(r,q)\rightarrow \cU(r,q)$ have
cardinality between $0$ and $d$.

We can pick a random element in $\cU(r,q)$ with uniform
distribution in the following way: we pick a random monic polynomial
 of
degree $r$ with coefficients in $\Fq$, with uniform
distribution. We  check whether  it is
irreducible. If it is, we output it. Otherwise we start again. This is
polynomial time in $r$ and $\log q$.

Given a random  element in $\cU(r,q)$ with uniform distribution,
we can compute the fiber of $x : \cP(r,q)\rightarrow \cU(r,q)$
above it and, provided  this fiber
 is non-empty,  pick a random element in it
with uniform distribution. If the
fiber is empty, we pick another element in $\cU(r,q)$ until we
find a non-empty fiber. At least one in every  $d\times(0.99)^{-1}$ fibers
is non-empty.
 We
thus define  a distribution $\mu$ on $\cP(r,q)$ and prove the following.

\begin{lemma}[A very rough measure]\label{lemma:random}
There is a unique measure $\mu$ on $\cP(r,q)$ such that
all non-empty fibers of the map $x : \cP(r,q)\rightarrow \cU(r,q)$
have the same measure, and all points in a given fiber
have the same measure.
There exists a probabilistic algorithm that picks a random element
in $\cP(r,q)$ with distribution $\mu$ in time polynomial in $d$
and $\log q$. For every subset $Z$ of $\cP(r,q)$ the measure
$\mu (Z)$ is related to  the uniform measure
$\frac{\#Z}{\#\cP(r,q)}$ by

$$\frac{\#Z}{d\#\cP(r,q)}\le   \mu (Z) \le \frac{d\#Z}{\#\cP(r,q)}.$$
\end{lemma}

Now let $\cD(r,q)$ be the set  of effective $\Fq$-divisors with  degree $r$ 
on $\cX$.
Since we have assumed  $q\ge 4g^2$ we know that   $\cX$ has at least one
$\Fq$-rational point. 
Let  $\Omega$ be  a degree $r$ effective divisor on
$\cX/\Fq$. We  associate to every $\alpha$ in $\cD(r,q)$ the class of
$\alpha-\Omega$ in $\cJ (\Fq)$. This defines a surjection $J_r :
\cD(r,q)\rightarrow \cJ (\Fq)$ with all
its fibers having  cardinality $\#\PP^{r-g}(\Fq)$.

So the set  $\cD( r,q)$  has cardinality
$\frac{q^{r-g+1}-1}{q-1}\#\cJ(\Fq)$. 

 So

$$\#\cP(r,q)  \le   \#\cD(r, q)  \le 
q^{r-g}\frac{1-\frac{1}{q^{r-g+1}}}{1-\frac{1}{q}}q^{g}(1+\frac{1}{\sqrt
q})^{2g}.$$

Since   $q\ge 4g^2$ we have   $\# \cD(r, q)\le 2eq^{r}$.

Assume $G$ is a finite group
and $\psi$  an epimorphism of groups $\psi : \cJ (\Fq)\rightarrow
G$.  We look for some divisor  $\Delta\in \cD(r,q)$ such that
$\psi(J_r(\Delta))\not= 0\in G$. Since all the fibers of $\psi\circ J_r$
have the same cardinality, the fiber above $0$ has at most
$\frac{2eq^{r}}{\#G}$ elements.
So the number of prime divisors $\Delta\in \cP(r,q)$ such that
$\psi(J_r(\Delta))$ is not $0$ is at least
$q^{r}(\frac{0.99}{r}-\frac{2e}{\# G})$. We assume $\#G$ is
at least $12r$. Then at least half of the divisors in
$\cP(r,q)$ are not mapped onto $0$ by $\psi\circ J_r$.
The $\mu$-measure of the subset consisting of these elements
is at least $\frac{1}{2d}$.

So if we pick a random $\Delta$ in $\cP(r,q)$ with $\mu$-measure as in
lemma \ref{lemma:random}, the probability of 
success is at least $\frac{1}{2d}$. If we make $2d$ trials, the probability
of success is $\ge 1-\exp(-1)\ge \frac{1}{2}$.

\begin{lemma}[Finding non-zero classes]\label{lemma:montecarlo}
There exists a probabilistic (Monte-Carlo) algorithm that takes as
input

\begin{enumerate}
\item   a degree $d$ and geometric genus $g$  plane projective 
absolutely irreducible reduced curve $C$ over
$\Fq$, such that $q\ge 4g^2$,

\item the smooth model $\cX$ of $C$,

\item a degree $g$ effective divisor $\omega$, as origin, 

\item an epimorphism   $\psi : \Pic^0(\cX/\Fq)\rightarrow G$  (that need not be computable) such that the cardinality of
$G$ is at least $\max(48g,24d,720)$,
\end{enumerate}
\noindent 
and outputs  a  sequence of $2d$
elements in $\Pic^0(\cX/\Fq)$  such that at least one of them is not in
the kernel of $\psi$ with probability $\ge \frac{1}{2}$. The algorithm
is polynomial time in $d$ and $\log q$.
\end{lemma}

As a special case we take $G=G_0=\cJ (\Fq)$ and $\psi=\psi_0$ the identity.
Applying lemma \ref{lemma:montecarlo} we find  a sequence of elements
in $\cJ (\Fq)$ out of which one at least is non-zero (with high
probability). We take  $G_1$
to be quotient of $G$ by the subgroup generated by these elements and
$\psi_1$ the quotient map. Applying the lemma again we construct
another sequence of elements in $\cJ (\Fq)$ out of which one at
least is not in $G_0$ (with high probability). We go on like that and produce a sequence of
subgroups in $\cJ (\Fq)$
that increase with constant probability until the index in
$\cJ (\Fq)$ becomes smaller than $\max(48g,24d,720)$. Note that
every step in this method is probabilistic: it succeeds with
some probability, that can be made very high (exponentially close
to $1$) while keeping a polynomial overall complexity.

\begin{lemma}[Finding an almost  generating set]\label{lemma:generators}
There exists a probabilistic (Monte-Carlo) algorithm that takes as
input

\begin{enumerate}
\item   a degree $d$ and geometric genus $g$  plane projective 
absolutely irreducible  reduced curve $C$ over
$\Fq$, such that $q\ge 4g^2$,

\item the smooth model $\cX$ of $C$,

\item a degree $g$ effective divisor $\omega$, as origin, 
\end{enumerate}
\noindent
and outputs  a  sequence of 
elements in $\Pic^0(\cX/\Fq)$   that generate a subgroup of index at most 

$$\max(48g,24d,720)$$
\noindent 
with probability $\ge \frac{1}{2}$. The algorithm
is polynomial time in $d$ and $\log q$.
\end{lemma}

Note that we do not catch the whole group $\cJ(\Fq)$ of rational points but
a subgroup $\cA$ with  index at most
$\iota = \max ( 48g, 24d, 720)$. This is a small
but annoying gap. In the sequel we shall try
to compute the $\ell$-torsion of the group $\cJ(\Fq)$
of rational points.
Because of the small gap in the above lemma, we may miss
some $\ell$-torsion points if $\ell$ is smaller than 
$\iota$. However, let $k$ be an integer
such that $\ell^k > \iota$. And let $x$ be a point
of order $\ell$
in $\cJ(\Fq)$. Assume there exists 
a point  $y$ in $\cJ(\Fq)$ such that
$x=\ell^{k-1}y$. The group $<y>$ generated by $y$ and the group
$\cA$ have non-trivial intersection because the 
product of their orders is bigger than the order of
$\cJ(\Fq)$. Therefore  $x$ belongs to $\cA$.

Our strategy for computing $\cJ(\Fq)[\ell]$ will be to find
a minimal field extension $\FQ$ of $\Fq$ such that all points
in $\cJ(\Fq)[\ell]$ are divisible by $\ell^{k-1}$ in $\cJ(\FQ)$. We then
shall apply the above lemma to $\cJ(\FQ)$. To finish
with,  we shall have to compute  $\cJ(\Fq)$ as a subgroup
of $\cJ(\FQ)$. To this end, we shall use the Weil pairing.

\section{Pairings}\label{section:pairings}

Let $n$ be a prime to $p$ integer and $\cJ$ a jacobian variety over $\Fq$. 
 The Weil
pairing relates the full $n$-torsion subgroup $\cJ(\Fqb)[n]$ with
itself. It can be defined using Kummer theory and is geometric in
nature.
The Tate-Lichtenbaum-Frey-R{\"u}ck pairing  is more cohomological and relates
the $n$-torsion $\cJ(\Fq)[n]$ in the group of $\Fq$-rational points
and the quotient $\cJ(\Fq)/n\cJ(\Fq)$. In this section, we quickly
review the definitions and algorithmic properties of these pairings,
following
work by Weil, Lang, Menezes,
Okamoto, Vanstone, Frey and R{\"u}ck.

\medskip 

We first recall the definition of
Weil pairing following \cite{langab}. Let
$k$ be an algebraically  closed field with characteristic $p$.
For every abelian variety $A$ over $k$, we denote by $Z_0(A)_0$
the group of $0$-cycles with degree $0$ and by $S : Z_0(A)_0\rightarrow A$ 
the summation map, that associates to every $0$-cycle of degree $0$
the corresponding sum in $A$.

Let $V$ and $W$ be two projective non-singular  irreducible and
reduced
varieties over $k$, and let 
$\alpha : V \rightarrow A$ and $\beta : W \rightarrow B$ be the
 canonical
maps into their Albanese varieties. Let $D$ be a correspondence on
$V\times W$.
Let $n\ge 2$ be a prime to $p$ integer. Let
$\agot$ (resp. $\bgot$) be a $0$-cycle of degree $0$ on $V$
(resp. $W$) and let $a=S(\alpha(\agot))$ (resp. $b=S(\beta(\bgot))$) 
be the associated
point in $A$ (resp. $B$). Assume $na=nb=0$. 
The  Weil
pairing $e_{n,D}(a,b)$ is defined 
in  \cite[VI, \S 4, Theorem 10]{langab}.
It is an $n$-th root of unity in $k$. It
depends  linearly in
$a$, $b$ and $D$.

Assume $V=W=\cX$ is a smooth projective  irreducible and 
reduced  curve over $k$ and $A=B=\cJ$ is its
jacobian and $\alpha=\beta=f :\cX \rightarrow \cJ$
is the Jacobi map (once  an origin on $\cX$ has been chosen). If we take $D$ to
be the diagonal on $\cX\times \cX$ we define a pairing $e_{n,D}(a,b)$
that
will be denoted $e_n(a,b)$ or $e_{n,\cX}(a,b)$. It does not depend on
the origin for the Jacobi map. It is non-degenerate.

The  jacobian $\cJ$ is principally polarized. We have
an isomorphism $\lambda : \cJ \rightarrow \hat \cJ$
between $\cJ$ and its dual $\hat \cJ$. If
$\alpha$ is an  endomorphism $\alpha : \cJ\rightarrow \cJ$,
we denote by ${}^t\alpha$
 its transpose  ${}^t\alpha : \hat \cJ
\rightarrow \hat \cJ$. If $D$ is a divisor on $\cJ$ that
is algebraically equivalent to zero, the image by ${}^t \alpha$
of the linear equivalence class of $D$ is the linear equivalence
class of the inverse image $\alpha^{-1}(D)$. See
\cite[V, \S 1]{langab}. The Rosati dual of $\alpha$ is defined to
be $\alpha^*=\lambda ^{-1}\circ {}^t\alpha \circ \lambda$. The 
map $\alpha \rightarrow \alpha^*$ is an involution, and 
$\alpha^*$ is the adjoint of $\alpha$ for the Weil pairing

\begin{equation}\label{equation:adjoint}
e_{n,\cX}(a,\alpha(b))=e_{n,\cX}(\alpha^*(a),b)
\end{equation}
\noindent according to \cite[VII, \S 2, Proposition 6]{langab}. 

If $\cY$ is another smooth projective irreducible and
reduced curve over $k$ and $\cK$ its jacobian and
$\phi : \cX\rightarrow \cY$ a non-constant map with degree $d$, and
$\phi^* : \cK\rightarrow \cJ$ the associated map between jacobians,  then
for $a$ and $b$ of order dividing $n$ in $\cK$ one has
$e_{n,\cX}(\phi^*(a),\phi^*(b))=e_{n,\cY}(a,b)^d$.

\medskip 

The Frey-R{\"u}ck pairing can be constructed from the Lichtenbaum version
of Tate's
pairing \cite{lichtenbaum} as was shown in \cite{freyruck}. Let $q$ be a
power of $p$. 
Let  again $n\ge 2$ be an   integer  prime to $p$ 
and $\cX$  a smooth projective absolutely
irreducible reduced curve over  $\Fq$. Let $g$ be the genus of $\cX$.
We assume $n$ divides $q-1$. 
Let $\cJ$ be the jacobian of $\cX$. The Frey-R{\"u}ck
pairing $\{,\}_n: \cJ (\Fq)[n]\times \cJ (\Fq)/n\cJ (\Fq)\rightarrow
\Fqs/(\Fqs)^n$
is defined as follows.  We take a class of order dividing $n$ in
$\cJ (\Fq)$. Such a class can be represented by an
$\Fq$-divisor $D$ with degree $0$. We take  a class in $\cJ (\Fq)$ and pick a degree zero
$\Fq$-divisor $E$  in this class, that we assume to be disjoint to $D$. The
pairing evaluated at the classes $[D]$ and $[E] \bmod n$ is $\{[D],[E]\bmod n \}_n=f(E)\bmod {(\Fqs)^n}$ where $f$
is any function with divisor $nD$.
This is a non-degenerate pairing. 

\medskip 

We now explain how one can compute the Weil pairing, following work by
Menezes, Okamo\-to, Vanstone, Frey and R{\"u}ck. The
Tate-Lichtenbaum-Frey-R{\"u}ck
pairing can be computed similarly.

As usual, we assume we are given a
degree $d$
plane model $C$ for $\cX$. 
Assume $\agot$ and $\bgot$ have
disjoint
support (otherwise we may replace $\agot$ by some linearly
equivalent
divisor using the explicit moving lemma \ref{lemma:moving}.)
We compute  a function $\phi$  with divisor $n\agot$. 
We similarly compute  a function $\psi$  with divisor $n\bgot$.
Then $e_n(a,b)=\frac{\psi(\agot)}{\phi(\bgot)}$. This algorithm is
polynomial in the degree $d$  of $C$ and the order $n$ of the divisors,
provided
the initial divisors $\agot$ and $\bgot$ are given as differences between effective
divisors with  polynomial degree  in $d$.

Using an idea that appears in a paper by Menezes, Okamoto and Vanstone
\cite{MOV} in the context of elliptic curves, and in \cite{freyruck}
for general curves, one can make this algorithm polynomial in $\log n$ in the
following way.
We write $\agot = \agot_0=\agot_0^+-\agot_0^-$ where $\agot_0^+$ and $\agot_0^-$ are effective
divisors. Let $\phi$  be the function computed in the above simple minded
algorithm. One has $(\phi)= n\agot_0^+-n\agot_0^-$. We want to express $\phi$ as a
product of small degree functions. We use a variant of fast
exponentiation. Using lemma \ref{lemma:moving} we compute a divisor 
$\agot_1= \agot_1^+-\agot_1^-$ and a function
$\phi_1$ such that $\agot_1$ is disjoint to $\bgot$ and
$(\phi_1)=\agot_1-2\agot_0$ and such that the degrees of $\agot_1^+$
and 
$\agot_1^-$ are
$\le 6gd (\log_q(\deg (\bgot))+1)$.
We go on  and compute, for $k\ge 1$ an integer, a divisor 
$\agot_k= \agot_k^+-\agot_k^-$ and a function
$\phi_k$ such that $\agot_k$ is disjoint to $\bgot$ and
$(\phi_k)=\agot_k-2\agot_{k-1}$ and such that the degrees of $\agot_k^+$
and 
$\agot_k^-$ are
$\le 6gd ( \log_q(\deg (\bgot)) +1)$.
We write the  base $2$ expansion of $n=\sum_i\epsilon_k2^k$  with
$\epsilon_k\in \{0,1\}$.
We compute the function $\Psi$ with divisor $\sum_k \epsilon_k\agot_k$.
We claim that the function $\phi$ can be written as a product of the
$\phi_k$, for $k\le \log_2n$,  and $\Psi$  with suitable integer exponents
bounded by $ n$ in absolute value. Indeed we write $\mu_1=\phi_1$, $\mu_2=\phi_2\phi_1^2$,
$\mu_3=\phi_3\phi_2^2\phi_1^4$ and so on. We have $(\mu_k)=\agot_k-2^k\agot$ and
${\Psi^{}}{\prod_k\mu_k^{-\epsilon_k}}$ has divisor $n\agot$ so is
the $\phi$ we were looking for.

\begin{lemma}[Computing the Weil pairing]\label{lemma:weilp}
There exists an algorithm that on input an   integer $n\ge
2$   prime to $q$ and 
a degree $d$  absolutely irreducible  reduced plane  projective curve $C$
over $\Fq$  and its  smooth
model $\cX$ 
and  two  $\Fq$-divisors on
 $\cX$,   denoted $\agot=\agot^+-\agot^-$ and $\bgot=\bgot^+-\bgot^-$, with
degree $0$,  and  order dividing $n$ in the jacobian,  computes the
Weil
pairing $e_n(\agot,\bgot)$ in time  polynomial 
in  $d$,  $\log q$,  $\log n$ and the degrees of $\agot^+$, $\agot^-$,
 $\bgot^+$, $\bgot^-$, the positive and negative parts of $\agot$ and
$\bgot$.
\end{lemma}

\begin{lemma}[Computation of  Tate-Lichtenbaum-Frey-R{\"u}ck
    pairings]\label{lemma:tatep}
There exists an algorithm that on input an integer $n\ge
2$  dividing $q-1$ and 
a degree $d$  absolutely irreducible reduced  plane  projective curve $C$
over $\Fq$  and its  smooth
model $\cX$ 
and  two $\Fq$-divisors on
 $\cX$,   denoted $\agot=\agot^+-\agot^-$ and $\bgot=\bgot^+-\bgot^-$, with
degree $0$,  and   such that the class of $\agot$ has order dividing $n\ge 2$ in the jacobian,  computes the
Tate-Lichtenbaum-Frey-R{\"u}ck 
pairing $\{\agot,\bgot\}_n$ in time polynomial
in  $d$,  $\log q$,  $\log n$ and the degrees of $\agot^+$, $\agot^-$,
 $\bgot^+$, $\bgot^-$, the positive and negative parts of $\agot$ and
$\bgot$.
\end{lemma}

\section{Divisible groups}\label{section:divisiblegroups}

Let $\Fq$ be a finite field with characteristic $p$ and
let  $\cX$ be a projective smooth absolutely irreducible  reduced algebraic
curve over $\Fq$. Let $g$ be the genus of $\cX$ and let $\ell\not = p$ be a
prime integer.  We assume $g\ge 1$.
Let $\cJ$ be the jacobian
of $\cX$ and let  $\End(\cJ/\Fq)$ be the ring of
endomorphisms of $\cJ$ over $\Fq$. Let $F_q$
be the Frobenius endomorphism.
In this section we study the action of $F_q$ on $\ell^k$-torsion
points of $\cJ$. We first consider the whole $\ell^k$-torsion
group. We then restrict to some well chosen subgroups
where this action is more amenable.

Let  $\chi(X)$  be the characteristic polynomial  of 
$F_q \in  \End(\cJ/\Fq)$.
The Rosati dual  to $F_q$  is
$q/F_q$. 
Let $\cO = \ZZ[X]/\chi(X)$ and
$\cO_\ell=\ZZ_\ell[X]/\chi(X)$.  We set $\varphi_q=X\bmod
\chi(X)\in \cO$.
Mapping  $\varphi_q$ onto $F_q$ defines an epimorphism from the ring
$\cO$  onto $\ZZ[F_q]$. 
In order to control the degree of the field of 
definition of $\ell^k$-torsion points we shall bound the order
of $\varphi_q$ in $(\cO/\ell^k\cO)^*$.

We set  $\cU_1=(\cO/\ell\cO)^*=(\Fl[X]/\chi(X))^*$.
Let the prime
factorization of $\chi(X)\bmod \ell$ be  $\prod_i\chi_i(X)^{e_i}$ with
$\deg(\chi_i)=f_i$. The order of $\cU_1$ is
$\prod_i\ell^{(e_i-1)f_i}(\ell^{f_i}-1)$. Let $\gamma$ be the
smallest integer such that  $\ell^{\gamma}$ is bigger than or equal
to  $2g$. Then the exponent of the group $\cU_1$ divides 
$A_1=\ell^{\gamma}\prod_i(\ell^{f_i}-1)$.
We set
$B_1=\prod_i(\ell^{f_i}-1)$
and $C_1=\ell^{\gamma}$. There is a unique polynomial $M_1(X)\in \ZZ[X] $
with degree $< 2g$  such that $\frac{\vp_q^{A_1}-1}{\ell}=M_1(\vp_q)\in
\cO$.

Now  for every positive 
integer $k$,  the element  $\vp_q$ belongs to the unit
group $\cU_k=(\cO/\ell^k\cO)^*$   of the quotient algebra
$\cO/\ell^k\cO=\ZZ[X]/(\ell^k,\chi(X))$.
The prime
factorization of $\chi(X)\bmod  \ell$ is lifted modulo $\ell^k$ as
$\prod_i\Xi_i(X)$ with $\Xi_i$ monic and 
$\deg(\Xi_i)=e_if_i$,  and  the order of $\cU_k$ is
$\prod_i\ell^{f_i(ke_i-1)}(\ell^{f_i}-1)$. The exponent of the latter
group divides
$A_k=A_1\ell^{k-1}$. So we 
set $B_k=B_1=\prod_i(\ell^{f_i}-1)$
and $C_k=C_1\ell^{k-1}=\ell^{k-1+\gamma}$. There is a unique polynomial
$M_k(X)
\in \ZZ[X] $
with degree $<  \deg(\chi)$  such that $\frac{\vp_q^{A_k}-1}{\ell^k}=M_k(\vp_q)\in
\cO$. 

For every integer $N\ge 2$ we can compute $M_k(X)\bmod N$ from
$\chi(X)$ in probabilistic polynomial time
in $\log q$, $\log \ell$, $\log N$, $k$, $g$:  we first factor
$\chi(X)\bmod \ell$ then compute the $\chi_i$ and the $e_i$ and $f_i$. We compute $X^{A_k}$ modulo $(\chi(X), \ell^kN)$
using fast exponentiation. We remove $1$ and divide by $\ell^k$.

\begin{lemma}[Frobenius and  $\ell$-torsion]\label{lemma:frobtorsion}
 Let $k$ be a
positive integer and $\ell\not = p$ a prime. Let $\chi(X)$ be the
characteristic polynomial of the Frobenius $F_q$ of $\cJ/\Fq$.
Let $e_i$ and $f_i$ be the multiplicities and inertiae in the prime
decomposition
of $\chi(X)\bmod \ell$. Let $\gamma$ be the smallest integer such that
$\ell^\gamma$ is bigger than or equal to  $2g$. Let  
$B=\prod_i(\ell^{f_i}-1)$. Let 
 $C_k=\ell^{k-1+\gamma}$ and $A_k=BC_k$. 
The $\ell^k$-torsion in $\cJ$ splits completely over the degree $A_k$
extension of $\Fq$. There is a degree
$<  2g$ polynomial $M_k(X)\in \ZZ[X]$ such that
$F_q^{A_k}=1+\ell^kM_k(F_q)$.
For every integer $N$ one  can compute such a $M_k(X)\bmod N$ from
$\chi(X)$ in probabilistic 
polynomial time 
in $\log q$, $\log \ell$, $\log N$, $k$, $g$.
\end{lemma}

In order to state sharper results  it is convenient to introduce
$\ell$-divisible subgroups inside  the $\ell^\infty$-torsion of a
jacobian $\cJ$,
that may
or may not correspond to subvarieties. We now see how to define such
subgroups and control their rationality properties.

\begin{lemma}[Divisible group]\label{definition:divisible}
Let  $\Pi :
J[\ell^\infty]\rightarrow J[\ell^\infty]$ be  a group
homomorphism whose restriction to its image $\GG$  is a bijection. 
Multiplication by $\ell$ is then  a surjection from $\GG$ to itself.
We denote by $\GG[\ell^k]$ the $\ell^k$-torsion in $\GG$.  There is an integer $w$ such that $\GG[\ell^k]$ is a
free $\ZZ/\ell^k\ZZ$ module of rank $w$ for every $k$.
We assume that  $\Pi$ commutes with the Frobenius endomorphism  $F_q$.
We then say $\GG$ is {\rm the
divisible group associated with $\Pi$}. From Tate's theorem \cite{tate}  $\Pi$ is  induced by some 
endomorphism in $\End (\cJ/\Fq)\otimes_\ZZ\ZZ_\ell$ and  we can define
$\Pi^*$  the Rosati dual of 
$\Pi$  and denote by  $\GG^*=\IM(\Pi^*)$ the associated
divisible group,  that we  call  the adjoint  of $\GG$.
\end{lemma}

\begin{remark}
The dual $\GG^*$ does not only depend on $\GG$. It may depend
on $\Pi$ also. 
\end{remark}

\begin{remark}
We may equivalently 
 define $\Pi^*$ as the dual of $\Pi$ for the Weil pairing.
See formula (\ref{equation:adjoint}).
\end{remark}

We now give an example of divisible group.
Let $F(X)=F_1(X)$ and $G(X)=G_1(X)$ be two
monic coprime polynomials in $\Fl[X]$ such that $\chi (X)=F_1(X)G_1(X)\bmod \ell$. From
Bezout's theorem we have two polynomials $H_1(X)$ and $K_1(X)$ in
$\Fl[X]$  such that $F_1H_1+G_1K_1=1$ and $\deg(H_1)<\deg(G_1)$ and $\deg(K_1)<\deg(F_1)$. From Hensel's lemma, for every
positive integer $k$ there exist four polynomials $F_k(X)$, $G_k(X)$,
$H_k(X)$ and $K_k(X)$ in $(\ZZ/\ell^k\ZZ)[X]$ such that $F_k$ and
$G_k$ are monic and
$\chi(X)=F_k(X)G_k(X)\bmod{\ell^k}$ and
$F_kH_k+G_kK_k=1\bmod{\ell^k}$  and $\deg(H_k)<\deg(G_1)$ and
$\deg(K_k)<\deg(F_1)$ and  $F_1=F_k\bmod \ell$, $G_1=G_k\bmod \ell$,
$H_1=H_k\bmod \ell$, $K_1=K_k\bmod \ell$.
The sequences $(F_k)_k$, $(G_k)_k$, $(H_k)_k$, $(K_k)_k$  converge in
$\ZZ_\ell[X]$ to $F_0$, $G_0$, $H_0$, $K_0$. 

If we substitute $F_q$ for $X$ 
 in $F_0H_0$ we obtain  a map $\Pi_G : \cJ[\ell^\infty]\rightarrow
\cJ[\ell^\infty]$ and similarly, if we substitute  $F_q$  for
$X$ in
$G_0K_0$ we obtain  a map $\Pi_F$.  It is clear that $\Pi_F^2=\Pi_F$
and $\Pi_G^2=\Pi_G$ and $\Pi_F+\Pi_G=1$ and $\Pi_F\Pi_G=0$. We call
$\GG_F=\IM(\Pi_F)$ and $\GG_G=\IM(\Pi_G)$
the associated supplementary $\ell$-divisible groups.

\begin{definition}[Characteristic subspaces]\label{definition:characteristic}
For every non-trivial monic factor $F(X)$ of $\chi(X)\bmod \ell$ such that the
cofactor 
$G=\chi/F\bmod\ell$ is prime to $F$, we write  $\chi=F_0G_0$ the
corresponding factorization in  $\ZZ_\ell[X]$. The $\ell$-divisible group
$\GG_F$ is
called the $F_0$-torsion in $\cJ[\ell^\infty]$ and is denoted 
$\cJ[\ell^\infty,F_0]$. It is the characteristic subspace of $F_q$
associated with the factor $F$. If $F=(X-1)^e$ is the largest  power
of $X-1$ dividing $\chi(X)\bmod \ell$ we abbreviate
$\GG_{(X-1)^e}=\GG_1$. If  $F=(X-q)^e$ then 
we write  similarly $\GG_{(X-q)^e}=\GG_{q}= \GG_1^*$.
\end{definition}

We now compute fields of definitions for torsion points inside
such divisible groups.
The action of $F_q$ on the $\ell^k$-torsion $\GG_F[\ell^k]=\cJ[\ell^k,
F_0]$ inside $\GG_F$
factors through the smaller ring $\cO_\ell/(\ell^k,F_0(\varphi_q))=\ZZ_\ell[X]/(\ell^k,F_0)$. We
deduce the following.

\begin{lemma}[Frobenius and  $F_0$-torsion]\label{lemma:frobtorsionGF}
 Let $k$ be a
positive integer and $\ell\not = p $ a prime. Let $\chi(X)$ be the
characteristic polynomial of the Frobenius $F_q$ of $\cJ$. Let $\chi =FG\bmod\ell$ with $F$ and $G$
monic coprime.
Let $e_i$ and $f_i$ be the multiplicities and inertiae in the prime
decomposition of $F(X)\bmod \ell$. Let $\gamma$ be the smallest integer such that
$\ell^\gamma$ is bigger than or equal to  $2g$. Let  
$B(F)=\prod_i(\ell^{f_i}-1)$. Let $C_k(F)=\ell^{k-1+\gamma}$ and $A_k(F)=B(F)C_k(F)$. 
The $\ell^k$-torsion in $\GG_F$ splits completely over the degree $A_k(F)$
extension of $\Fq$. There is a degree
$<  \deg(F)$ polynomial $M_k(X)\in \ZZ_\ell [X]$ such that
$\Pi_FF_q^{A_k(F)}=\Pi_F+\ell^k\Pi_FM_k(F_q)$.
For every power  $N$ of $\ell$,  one  can compute such an $M_k(X)$
modulo $N$ from
$\chi(X)$ and $F(X)$ in probabilistic 
polynomial time 
in $\log q$, $\log \ell$, $\log N$, $k$, $g$.
\end{lemma}

If we take for $F$ the largest power of $X-1$ dividing $\chi(X)\bmod \ell$  in
the above lemma, we can take   $B(F)=1$ so
$A_k(F)$ is an $\ell$ power $\le 2g\ell ^{k}$.

If we take for  $F$ the largest power of $X-q$ dividing $\chi(X)\bmod \ell$ in
the above lemma, we have
 $B(F)= \ell -1$ so
$A_k(F)$ is  $\le 2g(\ell-1)\ell ^{k}$.

So the characteristic  spaces associated with the eigenvalues $1$ and
$q$ split completely over small degree extensions of $\Fq$.

\section{The Kummer map}\label{section:kummer}

Let $\cX$ be a smooth projective absolutely  irreducible reduced  curve over
$\Fq$ of  genus $g$ and
$\cJ$ the jacobian of $\cX$.
Let $n\ge 2$ be an integer dividing $q-1$. We assume $g\ge 1$.
In this section, we construct a convenient surjection from
$\cJ(\Fq)$ to $\cJ(\Fq)[n]$.

If $P$ is in $\cJ(\Fq)$ we take some 
$R\in \cJ(\Fqb)$ such that $nR=P$ and form the $1$-cocycle
$({}^\sigma R -R)_\sigma$ in $H^1(\Fq, \cJ [n])$. Using the Weil pairing we
deduce an element 

$$\Box \mapsto (e_n({}^\sigma R -R,\Box))_\sigma$$
\noindent  in

$$\Hom
(\cJ [n](\Fq),H^1(\mmu_n))=\Hom(\cJ [n](\Fq),\Hom(\Gal(\Fq),\mmu_n)).$$

The map that sends $P\bmod n\cJ (\Fq)$ to $\Box\mapsto  (e_n({}^\sigma R -R,\Box))_\sigma$ is injective because the
Frey-R{\"u}ck pairing is non-degenerate. We observe that
$\Hom(\Gal(\Fq),\mmu_n)$
is isomorphic to $\mmu_n$: giving
an
 homomorphism from $\Gal(\Fq)$ to
$\mmu_n$ is equivalent to giving
the image of the Frobenius generator
$F_q$. We obtain a bijection $T_{n,q}$ from
$\cJ (\Fq)/n\cJ (\Fq)$
to the dual $\Hom(\cJ [n](\Fq),\mmu_n)$ of  $\cJ [n](\Fq)$ that we call the {\it Tate
map}. It maps  $P$ onto $\Box \mapsto e_n({}^{F_q} R -R,\Box)$. If
$\cJ[n]$ splits completely over $\Fq$ we set $K_{n,q}(P)={}^{F_q} R -R$
and define a bijection $K_{n,q} : \cJ(\Fq)/n\cJ(\Fq) \rightarrow
\cJ[n](\Fq)=\cJ[n]$ that we call the {\it Kummer map}.

\begin{definition}[The Kummer map] Let $\cJ/\Fq$ be a jacobian and $n\ge 2$ an
integer. Assume $\cJ[n]$ splits completely over $\Fq$. For $P$
in $\cJ(\Fq)$ we choose any $R$ in $\cJ(\Fqb)$ such that
$nR=P$ and we set $K_{n,q}(P)= {}^{F_q} R -R$.
This defines a bijection
$$K_{n,q} : \cJ(\Fq)/n\cJ(\Fq) \rightarrow
\cJ[n](\Fq)=\cJ[n].$$
\end{definition}

We now assume that $n=\ell^k$ is a  power of some prime
integer $\ell\not = p $. We also  make the (strong~!) assumption that $\cJ [n]$
splits completely over $\Fq$. 
We want to compute the Kummer map $K_{n,q}$
explicitly. Let $P$ be an
$\FF_{q}$-rational
point in $\cJ$. Let $R$ be such that $nR=P$. 
Since $F_q-1$ kills $\cJ [n]$, 
there is an $\Fq$-endomorphism $\kappa$
of $\cJ$ such that $F_q-1=n\kappa$.  We note that $\kappa$ belongs to 
$\ZZ[F_q]\otimes_\ZZ\QQ=\QQ[F_q]$  and therefore commutes with $F_q$. We have $\kappa
(P) = (F_q-1)(R)=K_{n,q}(P)$ and $\kappa (P)$ is $\FF_{q}$-rational.

Computing the Kummer map will be seen to be 
 very useful but it requires that
$\cJ [n]$ splits completely
over $\Fq$. In general, we shall have to base change to some extension
of $\Fq$.

Let $\chi(X)$ be the characteristic polynomial of $F_q$ and
let $B=\prod_i(\ell^{f_i}-1)$ where the $f_i$ are the degrees of prime
divisors of $\chi(X)\pmod \ell$. Let $\ell^\gamma$ be the smallest power of $\ell$ that is
bigger than or equal to $2g$.
Let $C_k=\ell^{\gamma+k-1}$ and $A_k=BC_k$.   Set $Q=q^{A_k}$.
  From lemma  \ref{lemma:frobtorsion} 
there is a polynomial $M_k(X)$ such that $F_Q=1+\ell^kM_k(F_q)$.
So, for $P$ an
$\FF_{Q}$-rational
point in $\cJ$ and  $R$ such that $nR=P$,  the Kummer map  $K_{n,Q}$
applied to $P$ is  $M_k(F_q)(P) = (F_Q-1)(R)=K_{n,Q}(P)$ and this is
an  $\FF_{Q}$-rational point.

\begin{lemma}[Computing the Kummer map]\label{lemma:kummer1}
Let $\cJ/\Fq$ be a jacobian. Let $g\ge 1$ be its dimension. 
Let  $\ell\not= p$ be a
prime integer and $n=\ell^k$ a power of $\ell$. 
Let $\chi(X)$ be the characteristic polynomial of $F_q$ and
let $B=\prod_i(\ell^{f_i}-1)$ where the $f_i$ are the degrees of prime
divisors of $\chi(X)\pmod \ell$. 
Let $\ell^\gamma$ be the smallest power of $\ell$ that is
bigger than or equal to $2g$.
Let $C_k=\ell^{\gamma+k-1}$ and $A_k=BC_k$.   Set $Q=q^{A_k}$ and  observe
that  $n$ divides $Q-1$ because $\cJ[n]$ splits completely over $\FF_Q$.
There exists an endomorphism $\kappa \in \ZZ[F_q]$ of $\cJ$ such that $n\kappa
=F_Q-1$
 and  for every $\FF_{Q}$-rational point $P$ and
any
$R$ with $nR=P$ one has  $\kappa
(P) = (F_{Q}-1)(R)=K_{n,Q}(P)$. This endomorphism $\kappa$ induces a bijection
between $\cJ (\FF_{Q})/n\cJ (\FF_{Q})$ and $\cJ [n](\FF_Q)=\cJ [n]$. Given $\chi(X)$
and a positive 
integer $N$ one can compute $\kappa \bmod N$  as a polynomial in $F_q$
with coefficients in $\ZZ/N\ZZ$ in probabilistic polynomial time in
$g$, $\log \ell$,  $\log q$, $k$, $\log N$.
\end{lemma}

This lemma is not of much use in practice because the 
 field  $\FQ$  is too big.
On the other hand, we may not be interested in the whole $n$-torsion in
$\cJ$ but just a small piece in it, namely the $n$-torsion
of a given divisible group.

So let $\ell\not =
p$ be a prime integer and $\GG$ an $\ell$-divisible group in
$\cJ[\ell^\infty]$
and $\Pi=\Pi^2 :
\cJ[\ell^\infty]\rightarrow \GG$ a projection onto it.
Let   $n=\ell^k$ and let $Q$ be a power of $q$ such that $\GG[n]$
splits completely over $\FF_Q$. 
Let $P$ be an
$\FQ$-rational
point in $\GG$. Let $R\in \GG(\Fqb)$ be such that $nR=P$.
We set $K_{\GG,n,Q} (P)= {}^{F_Q}R-R$ and define an isomorphism

$$K_{\GG,n,Q} : \GG(\FF_Q)/n\GG(\FF_Q) \rightarrow
\GG(\FF_Q)[n]=\GG[n].$$

In order to make this construction explicit, we now assume that there
exists some 
$\kappa \in \ZZ_\ell[F_q]$ such that 
$\Pi(F_Q-1-n\kappa )=0$.
Lemma \ref{lemma:frobtorsionGF} provides us with such a $Q$ and such a
$\kappa$ when 
$\GG=\cJ[\ell^\infty,F_0]$ is some
characteristic subspace.

We now can  compute this new  Kummer map $K_{\GG,n,Q}$. Let $P$ be an
$\FQ$-rational
point in $\GG$. Let $R\in \GG$ be such that $nR=P$. 
From $(F_Q-1-n\kappa )\Pi(R)=0=(F_Q-1-n\kappa )(R)$ we deduce that
$K_{\GG, n,Q}(P)=\kappa(P)$. Hence the 

\begin{lemma}[The Kummer map for a divisible group]\label{lemma:kummer2}
Let $\cJ/\Fq$ be a jacobian. Let $g$ be its dimension.
Let $\ell\not = p$ be a
prime integer and $n=\ell^k$ a power of $\ell$. We assume $g\ge 1$.
Let $\chi(X)$ be the characteristic polynomial of $F_q$. Assume
$\chi(X)=F(X)G(X) \bmod\ell$ with $F$ and $G$ monic coprime polynomials in
$\Fl[X]$ and let $\GG_F$ be the associated $\ell$-divisible group.
Let $B=(\ell-1)\prod_i(\ell^{f_i}-1)$ where the $f_i$ are the degrees of prime
divisors of $F(X)\pmod \ell$. 
 Let $\ell^\gamma$ be the smallest power of $\ell$ that is
bigger than or equal to $2g$.
Let $C_k=\ell^{k-1+\gamma}$ and $A_k=BC_k$. Set $Q=q^{A_k}$. From
lemma
\ref{lemma:frobtorsionGF} there exists an endomorphism  $\kappa \in \ZZ_\ell[F_q]$ such
that   $\Pi_F(n\kappa- F_Q+1)=0$
 and  for every $\FF_{Q}$-rational point $P\in \GG_F$ and
any
$R\in \GG_F$ with $nR=P$ one has  $\kappa
(P) = (F_{Q}-1)(R)=K_{\GG, n,Q}(P)$. This endomorphism $\kappa$ induces a bijection
between $\GG_F (\FF_{Q})/n\GG_F (\FF_{Q})$ and $\GG_F [n](\FQ)=\GG_F [n]$. Given $\chi(X)$
and $F(X)$ and  a power 
$N$ of $\ell$,  one can compute $\kappa \bmod N$  as a polynomial in $F_q$
with coefficients in $\ZZ/N\ZZ$ in probabilistic polynomial time in
$g$, $\log \ell$,  $\log q$, $k$, $\log N$.
\end{lemma}

\section{Linearization of torsion classes}\label{section:relations}

Let $C$ be a  degree $d$   plane projective absolutely irreducible
reduced 
curve $C$ over
$\Fq$ with geometric genus $g\ge 1$,
and assume  we are given  the smooth model $\cX$ of $C$. 
We also assume we are given  a degree $1$ divisor $O=O^+-O^-$ where $O^+$
and $O^-$ are effective, $\Fq$-rational and have degree bounded by an
absolute 
constant times $g$.

Let 
$\cJ$ be the  jacobian of $\cX$.
We  assume $\ell\not = p$ is a prime integer that divides
$\#\cJ(\Fq)$. Let $n=\ell^k$ be a power of $\ell$. We want
to describe $\cJ(\Fq)[\ell^k]$ by generators and relations.

If  $x_1$, $x_2$, \dots, $x_I$ are  elements in  a finite commutative group $G$ we let $\cR$
be the kernel of the map $\xi : \ZZ^I\rightarrow G$ defined by 
$\xi(a_1,\cdots,a_I)=\sum_ia_ix_i$. We call $\cR$ the {\it lattice of
relations} between the $x_i$.

We first give a very general and rough algorithm for computing
relations in any finite commutative group.

\begin{lemma}[Finding relations in  blackbox groups]\label{lemma:idiot}
Let $G$ be a finite and commutative group and let $x_1$, $x_2$, ...,
$x_I$ be elements in $G$. A basis for the lattice of relations between
the $x_i$ can be
computed at the expense of $3I\#G$ operations (or comparisons) in $G$.
\end{lemma}

We first compute and store all the multiples of $x_1$. So  we list $0$, $x_1$,
$2x_1$, \dots until we find the first multiple $e_1x_1$ that is equal to zero.
This gives us the relation $r_1=(e_1,0,\ldots,0)\in \cR$. This first
step requires at most $o=\#G$ operations in $G$  and $o$ comparisons.

We then compute successive multiples of $x_2$ until we find the first
one
$e_2x_2$ that is in $L_1=\{0,x_1,\ldots,(e_1-1)x_1\}$. This gives us a
second
relation $r_2$. The couple $(r_1,r_2)$ is a basis for the
 lattice of relations between $x_1$ and $x_2$. Using this lattice, we
 compute the list $L_2$ of elements in the group generated by $x_1$
 and $x_2$. This second
step requires at most $2o$ operations and $e_1e_2\le o$ comparisons.

We then compute successive multiples of $x_3$ until we find the first
one
$e_3x_3$ that is in $L_2$. This gives us a
third
relation $r_3$. The triple $(r_1,r_2,r_3)$ is a basis for the
 lattice of relations between $x_1$, $x_2$ and $x_3$. Using this lattice, we
 compute the list $L_3$ of elements in the group generated by $x_1$,
 $x_2$ and $x_3$. This third
step requires at most $2o$ operations and $o$ comparisons. And we go
on like this. \hfill
$\Box$

This is far from efficient unless the group is very small.

We  come back to the computation of generators and relations
for $\cJ(\Fq)[\ell^k]$. 

Let $B=\ell-1$. Let $\ell^\gamma$ be
the smallest power of $\ell$  that  is bigger than or equal to $2g$
and let $A_k= B\ell^{\gamma + k-1}$. We set $Q_k=q^{A_k}$. 

If we take  for $F$  a power of $X-1$ in definition \ref{definition:characteristic}
and lemma \ref{lemma:kummer2}  we obtain two surjective maps $\Pi_1 :
\cJ(\FF_{Q_k})[\ell^\infty]\rightarrow \GG_1(\FF_{Q_k})$ and
$K_{\GG_1,\ell^k,Q_k} : \GG_1(\FF_{Q_k}) \rightarrow \GG_1[\ell^k]$.

If we now take for  $F$ a power of $X-q$ in definition \ref{definition:characteristic}
and lemma \ref{lemma:kummer2}  we obtain two surjective maps $\Pi_{q} :
\cJ(\FF_{Q_k})[\ell^\infty]\rightarrow \GG_{q}(\FF_{Q_k})$ and
$K_{\GG_{q},\ell^k,Q_k} : \GG_{q}(\FF_{Q_k}) \rightarrow \GG_{q}[\ell^k]$.

There exists a unit $u$ in $\End (\cJ/\Fq)\otimes_\ZZ\ZZ_\ell$ such that
the Rosati dual $\Pi_1^*$ of $\Pi_1$  is 

$$\Pi_1^*=u\Pi_q.$$

Therefore  $\GG_q=\GG_1^*$ and
the restriction of the Weil pairing to $\GG_1[\ell^k]\times 
\GG_{q}[\ell^k]$  is non-degenerate.

If $Q_k\ge 4g^2$, we use  lemma \ref{lemma:generators} to
  produce a sequence
$\gamma_1$, \dots, $\gamma_I$ of
elements in  $\cJ(\FF_{Q_k})$ that generate (with high probability)
a subgroup of index at
most
$\iota = \max ( 48g, 24d, 720)$.
If $Q_k\le 4g^2$  we use lemma \ref{lemma:Stein} to produce
a sequence 
$\gamma_1$, \dots, $\gamma_I$ of
elements in  $\cJ(\FF_{Q_k})$ that generate it.

 Let $N$ be the largest divisor of $\#\cJ(\FF_{Q_k})$ which is prime to
$\ell$. 

We set $\alpha_i=K_{\GG_1,\ell^k,Q_k}(\Pi_1(N\gamma_i))$ and 
$\beta_i=K_{\GG_{q},\ell^k,Q_k}(\Pi_{q}(N\gamma_i))$.

The group $\cA_k$ generated by the $\alpha_i$ has index at most $\iota$
in
$\GG_1[\ell^k]$.
 The group $\cB_k$ generated by the $\beta_i$ has index at most $\iota$
in
$\GG_{q}[\ell^k]$.

Let $\ell^\delta$ be smallest power of $\ell$ that is bigger than
$\iota$ and assume $k>\delta$.  Then $\cA_k$ contains
$\GG_1[\ell^{k-\delta}]$.

We now explain how to compute the lattice of relations  between
given elements $\rho_1$, \dots , $\rho_J$ in $\GG_1[\ell^k]$. We denote by
$\cR$ this lattice. Recall
the restriction of the Weil pairing  to $\GG_1[\ell^k]\times 
\GG_{q}[\ell^k]$  is a non-degenerate pairing

$$e_{\ell^k} : \GG_1[\ell^k]\times 
\GG_{q}[\ell^k]\rightarrow \mmu_{\ell^k}.$$

We fix an isomorphism 
between the group $\mmu_{\ell^k}(\Fqb)=
\mmu_{\ell^k}(\FF_{Q_k})$ of $\ell^k$-th roots of unity and
$\ZZ/\ell^k\ZZ$. 
Having  chosen the preimage of $1\bmod \ell^k$, computing
this isomorphism is a problem called {\it discrete logarithm}.
We can  compute this discrete logarithm
 by exhaustive search at the expense of $O(\ell^k)$ operations in $\FF_{Q_k}$.
There exist more efficient algorithms, but we don't need them for our
complexity
estimates.

We regard  the matrix $(e_{\ell^k}(\beta_i,\rho_j))$ as a matrix with
$I$ rows, $J$ columns and coefficients in $\ZZ/\ell^k\ZZ$. This
matrix defines a morphism from $\ZZ^J$ to $(\ZZ/\ell^k\ZZ)^I$ whose
kernel
is a lattice $\cR'$ that contains $\cR$. The index of $\cR$ in $\cR'$
is at most $\iota$. Indeed $\cR'/\cR$ is isomorphic to the orthogonal
complement 
of $\cB_k$ in $<\rho_1,\ldots,\rho_J>\subset \GG_1[\ell^k]$. So it  has order $\le \iota$.
We then compute a basis of $\cR'$. This boils down to computing
the kernel of an $I\times (J+I)$ integer matrix with
entries bounded by $\ell^k$. This can be done
by putting this 
matrix in Hermite normal form (see \cite[2.4.3]{cohen}). The complexity is polynomial in 
$I$, $J$ and $k\log \ell$. See \cite{Havas},  \cite[2.4.3]{cohen} and \cite{vdk}.

Once given a basis of $\cR'$, the sublattice $\cR$ can be  computed
using lemma \ref{lemma:idiot} at the expense of $\le 3J\iota$ operations.

We apply this method to the generators $(\alpha_i)_i$ of $\cA_k$. Once
given the lattice $\cR$ of relations between the $\alpha_i$ it is a
matter of linear algebra to find a basis $(b_1,\dots,b_w)$  for
$\cA_k[\ell^{k-\delta}]=\GG_1[\ell^{k-\delta}]$.
The latter group is a rank $w$ free module over $\ZZ/\ell^{k-\delta}\ZZ$ and is acted
on by the $q$-Frobenius $F_q$. For every $b_j$ we can compute the
lattice
of relations between $F_q(b_j)$, $b_1$, $b_2$, \dots, $b_w$ and deduce
the matrix of $F_q$ with respect to
 the basis $(b_1,\dots,b_w)$. From this matrix 
we deduce
a nice generating set  for the kernel  of $F_q-1$ in $\GG_1[\ell^{k-\delta}]$. This
kernel is $\cJ[\ell^{k-\delta}](\Fq)$. 
We deduce the following.

\begin{theorem}\label{theorem:mainjacobi}
There is a probabilistic Monte-Carlo  algorithm that on  input 

\begin{enumerate}
\item   a degree $d$ and geometric genus $g$  plane projective 
absolutely irreducible reduced curve $C$ over
$\Fq$, 

\item the smooth model $\cX$ of $C$,

\item  a degree $1$ divisor $O=O^+-O^-$ where $O^+$
and $O^-$ are effective, $\Fq$-rational and have degree bounded by a
constant times $g$,

\item a prime $\ell$ different from the characteristic $p$ of $\Fq$
  and a power $n=\ell^k$ of $\ell$,
\item the zeta function of $\cX$;
\end{enumerate}
outputs 
a set $g_1$, \dots, $g_W$ of divisor classes in   the
Picard group of $\cX/\Fq$,   such that
the $\ell^k$ torsion $\Pic(\cX/\Fq)[\ell^k]$ is the direct product
of the $<g_i>$,  and the orders of the $g_i$ form a non-decreasing
sequence.
Every class $g_i$ is given by a divisor $G_i-gO$ in the class, where $G_i$ is
a degree $g$  effective $\Fq$-divisor on $\cX$.

The algorithm runs in probabilistic polynomial time in $d$, $g$,
$\log q$ and $\ell^k$. It outputs the correct answer with probability
$\ge \frac{1}{2}$. Otherwise, it may return either nothing or a strict
subgroup of $\Pic(\cX/\Fq)[\ell^k]$.

If one is given a degree zero $\Fq$-divisor $D=D^+-D^-$ of order
dividing $\ell^k$,
one can compute the coordinates of the class of $D$ in the basis
$(g_i)_{1\le i \le W}$ in polynomial time in $d$, $\log q$, $\ell^k$
and
the degree of $D^+$. These coordinates are integers $x_i$ such that $\sum_{1\le i\le W}x_ig_i=[D]$.
\end{theorem}

\section{An example: modular curves}\label{section:modularcurves}

In this section we consider a family of modular curves for which we
can
easily provide and study  a plane model. 
Let $\ell \ge 5$ be a prime. We set $d_\ell=\frac{\ell^2-1}{4}$ and
$m_\ell=\frac{\ell -1}{2}$. We denote by $\cX_\ell=X(2)_1(\ell)$ the
moduli of elliptic curves with full $2$-torsion plus  one non-trivial
$\ell$-torsion point. We first describe a homogeneous singular plane  model $C_\ell$ for
this curve. We enumerate the geometric points  on $\cX_\ell$ above every 
singularity of $C_\ell$ and  compute the conductor  $\Cgot_\ell$ using
the Tate elliptic  curve.

Let $\lambda$ be an indeterminate and form the Legendre elliptic curve
with equation $y^2=x(x-1)(x-\lambda)$. Call $\cT_\ell(\lambda,x)$ the
$\ell$-division  polynomial of this curve. It is 
 a polynomial in $\QQ[\lambda][x]$ with 
degree
$2d_\ell=\frac{\ell^2-1}{2}$ in $x$.

As a polynomial in $x$ we have

$$\cT_\ell(\lambda,x)=\sum_{0\le k \le
  2d_\ell}a_{2d_\ell-k}(\lambda)x^k$$
\noindent where  $a_0(\lambda)$ has degree $0$ in $\lambda$ so that
we  normalise  by setting
$a_{0}(\lambda)=\ell$. 


Let $\cF$ be a splitting field of $\cT_\ell(\lambda,x)$ over
$\QQ(\lambda)$. A suitable twist of the Legendre curve has
a point of order $\ell$ defined over $\cF$ (and the full two
torsion also). This proves that $\cF$ contains the function
field  $\QQ(\cX_\ell)$. Comparison of the degrees of
$\cF/\QQ(\lambda)$ and $\QQ(\cX_\ell)/\QQ(\lambda)$ shows that
the two fields $\cF$ and $\QQ(\cX_\ell)$ are equal and the polynomial
$\cT_\ell$ is irreducible in $\bar \QQ (\lambda)[x]$.

We can compute the $2d_\ell$ roots of $\cT_\ell(\lambda,x)$ in the field
$\Qb \{\{ \lambda^{-1} \}\}$ of
Puiseux series in $\lambda^{-1}$.
We set 

$$j=j(\lambda)=2^8\frac{(\lambda^2-\lambda+1)^3}{\lambda^2(\lambda
  -1)^2}=2^8\lambda^2(1-\lambda^{-1}+
3\lambda^{-2}+3\lambda^{-4}+\cdots)$$
\noindent  so that 
$j^{-1}=2^{-8}(\lambda^{-2}+\lambda^{-3}-2\lambda^{-4}-5\lambda^{-5}+\cdots)$.

We introduce Tate's  $q$-parameter,
defined implicitly by 

$$j=\frac{1}{q}+744+196884q+\cdots$$
\noindent  so that

\begin{eqnarray*}
q&=&j^{-1}+744j^{-2}+750420j^{-3}+\cdots\\
&=&  \frac{1}{256}\lambda^{-2} +
\frac{1}{256}\lambda^{-3} +
\frac{29}{8192}\lambda^{-4}+\frac{13}{4096}\lambda^{-5}+\cdots
\end{eqnarray*}

We set $x=x'+\frac{1+\lambda}{3}$ and $y'=y$ and  find the reduced Weierstrass equation for
the Legendre curve 

$$y'^2=x'^3-\frac{\lambda^2-\lambda+1}{3}x'-\frac{(\lambda-2)(\lambda+1)(2\lambda-1)}{27}.$$

We want to compare the latter curve and the Tate curve with equation

$$y''^2=x''^3-\frac{E_4(q)}{48}x''+\frac{E_6(q)}{864}$$
\noindent  where
$E_4(q)=1+240q+\cdots$ and $E_6(q)=1-504q+\cdots$. 

The quotient $\frac{E_4(q)(dq)^2}{(\lambda^2-\lambda+1)q^2}$ is a quadratic  differential on the
curve $X(2)$ with divisor $-2(0) -2(1)$ in the $\lambda$
coordinate. Examination of the leading terms of its expansion shows
that

$$E_4\left(\frac{dq}{q}\right)^2=\frac{4(\lambda^2-\lambda+1)(d\lambda)^2}{\lambda^2(1-\lambda)^2}$$
\noindent 
and similarly
$$E_6\left(\frac{dq}{q}\right)^3=\frac{4(\lambda-2)(\lambda+1)(2\lambda-1)(d\lambda)^3}{\lambda^3(1-\lambda)^3}.$$

We deduce the isomorphism $x'=\gamma ^2x''$  and $y'=\gamma^3y''$ with

$$\gamma^2=2\lambda(\lambda-1)\left(
\frac{dq}{qd\lambda}\right)=-4\lambda+2+\frac{3}{8}\lambda^{-1}+\frac{3}{16}\lambda
^{-2} + \cdots.$$

Set  $\zeta_\ell=\exp(\frac{2i\pi}{\ell})$.
For $a$ and $b$ integers such that either $b=0$ and
$1\le a \le \frac{\ell -1}{2}$ or $1\le  b \le  \frac{\ell-1}{2}$ and
$0\le a \le \ell -1$  we set
$w=\zeta_\ell^aq^{\frac{b}{\ell}}$ in the  expansion 

$$x''(w, q)=\frac{1}{12}+\sum_{n\in \ZZ}\frac{wq^n}{(1-wq^n)^2}
-2\sum_{n\ge 1}\frac{nq^n}{1-q^n}$$
\noindent and find

$$x''_{a,b}=\frac{1}{12}+\zeta_\ell^aq^{\frac{b}{\ell}}+O(q^{\frac{b+1}{\ell}})$$
\noindent if $b\not =0$,  and
$x''_{a,0}=\frac{1}{12}+ \frac{\zeta_\ell^a}{(1-\zeta_\ell^a)^2}+O(q)$.

So

$$x_{a,b}=\gamma^2
x''+\frac{1+\lambda}{3}=-4\zeta_\ell^a2^{\frac{-8b}{\ell}}\lambda^{1-\frac{2b}{\ell}}+O(\lambda^{1-\frac{2b+1}{\ell}})$$
\noindent if $b\not =0$ and $x_{a,0}= \frac{-4\zeta_\ell^a}{(1-\zeta_\ell^a)^2}\lambda+O(1)$.

The $x_{a,b}$ are the roots of $\cT_\ell(\lambda,x)$ in the field
$\Qb \{\{ \lambda^{-1} \}\}$ of Puiseux series.

We deduce that for  $1\le k\le \frac{\ell-1}{2}$ the  polynomial
$a_k(\lambda)$ has degree at most $k$. Further $a_{\frac{\ell
    -1}{2}}(\lambda)=2^{{\ell - 1}}(-\lambda)^{\frac{\ell -
    1}{2}}+O(\lambda^{\frac{\ell-3}{2}})$. For $k>\frac{\ell-1}{2}$
 the  polynomial
$a_k(\lambda)$ has degree $<k$ and $\le d_\ell$.

The coefficients in all the series expansions above are in
 $\ZZ [\frac{1}{6\ell}, \zeta_\ell, 2^{\frac{1}{\ell}}]$. The coefficients of  $\cT_\ell(\lambda, x)$
 are
in  $\ZZ[\frac{1}{6\ell}]$. In fact $\cT_\ell(\lambda,x)$ is in
$\ZZ[\lambda,x]$ but this is not needed here.

Since $\cT_\ell\in \QQ [\lambda, x]$ is absolutely irreducible, the
equation $\cT_\ell(\lambda,x)=0$  defines a plane absolutely irreducible
affine curve $\cC_\ell$. 
Let  $C_\ell\subset \PP^2$ be the 
projective plane
curve made of the zeroes of the  homogeneous polynomial
$\cT_\ell(\frac{\Lambda}{Y},\frac{X}{Y})Y^{2d_\ell}$.

For every geometric point $P$ on $\cX_\ell$ such that $\lambda (P) \not \in
\{0,1,\infty \}$,  the function $\lambda -\lambda(P)$ is a
uniformizing parameter at $P$. 
Further $x(P)$ is finite and $P$ is the only geometric point on  $\cX_\ell$ above
the point $(\lambda(P), x(P))$ of $\cC_\ell$. So the only possible
singularities of $C_\ell$ lie on one of the three lines with equations
$\Lambda =0$, $Y=0$ and $\Lambda-Y=0$.

The points at infinity are given by the degree $2d_\ell$ form
$$2^{\ell-1}(-1)^{\frac{\ell -1}{2}}\Lambda^{\frac{\ell
    -1}{2}}X^{\frac{\ell^2 -\ell}{2}}+\dots + \ell X^{\frac{\ell^2
    -1}{2}}=X^{\frac{\ell^2 -\ell}{2}}\prod_{0\le a\le \frac{\ell
    -1}{2}}(-4\Lambda-(\zeta_\ell^a+\zeta_\ell^{-a}-2)X).$$

We call $\Sigma_\infty = [1,0,0]$ the unique singular point  at
infinity and for every $1\le b\le \frac{\ell -1}{2}$ we call
$\sigma_{\infty, b}$ the point above $\Sigma_\infty$ on $\cX_\ell$
associated with  the orbit

$$\{x_{0,b}, x_{1,b}, \cdots, x_{\ell -1,b}
\}$$
\noindent  for the local monodromy group. We call $\mu_{\infty, a}$ the point on
$\cX_\ell$
corresponding to the expansion $x_{a,0}$.  
The ramification index of the covering map
$\lambda : \cX_\ell \rightarrow X(2)$ is $\ell$ at $\sigma_{\infty,b}$ and $1$
at $\mu_{\infty, a}$. Since $\ell-2b$ and $\ell$
are coprime, there exist two integers $\alpha_b$ and $\beta_b$ such
that $\alpha_b(\ell-2b)-\beta_b\ell=1$ and $1\le \alpha_b\le \ell-1$
and  $1\le \beta_b \le  \ell-1$. The monomial
$x^{\alpha_b}\lambda^{-\beta_b}\in \bQ(\cX_\ell)$
is a local parameter at $\sigma_{\infty,b}$. Of course,
$\lambda^{-\frac{1}{\ell}}$ is also a local parameter at this point,
and it is much more convenient, although it is not in $\bQ(\cX_\ell)$.

The morphism $\phi : \cX_\ell\rightarrow X_1(\ell)$ corresponding to forgetting
the $2$-torsion structure is Galois with group $\cS_3$ generated by
the two transpositions $\tau_{(0,\infty)}$ and $\tau_{(0,1)}$ defined in
homogeneous coordinates by

$$\tau_{(0,\infty)} : [\Lambda,X,Y]\rightarrow [Y,X,\Lambda]$$
\noindent   and 

$$\tau_{(0,1)}
: [\Lambda,X,Y]\rightarrow [Y-\Lambda,Y-X,Y].$$

We observe that these
act on  $\cX_\ell$, $\PP^2$ and $C_\ell$ in a way compatible
with the maps $\cX_\ell\rightarrow C_\ell$ and $C_\ell\subset \PP^2$.
We set  $\Sigma_0=\tau_{(0,\infty)}(\Sigma_\infty) = [0,0,1]$ and
$\Sigma_1=\tau_{(0,1)}(\Sigma_0) = [1,1,1]$.
We set $\sigma_{0, b}=\tau_{(0,\infty)}(\sigma_{\infty, b})$  and
$\sigma_{1, b}=\tau_{(0,1)}(\sigma_{0, b})$,  $\mu_{0, a}=\tau_{(0,\infty)}(\mu_{\infty, a})$  and $\mu_{1, a}=\tau_{(0,1)}(\mu_{0, a})$.

The  genus of $\cX_\ell$  is $g_\ell=\frac{(\ell - 3)^2}{4}=(m_\ell-1)^2$. The arithmetic
genus of $C_\ell$ is $g_a=(m_\ell^2+m_\ell-1)(2m_\ell^2+2m_\ell-1)$. We now compute the
conductor of $C_\ell$. Locally at $\Sigma_\infty$ the curve $C_\ell$
consists of $m_\ell$ branches (one for each point $\sigma_{\infty,b}$) that are cusps with equations

$$\left(\frac{X}{\Lambda}\right)^\ell=
-2^{2\ell-8b}\left(\frac{Y}{\Lambda}\right)^{2b}+\cdots$$

The
conductor of this latter 
cusp is $\sigma_{\infty,b}$ times $(\ell-1)(2b-1)$ which 
is  the next integer to the last gap
of the additive semigroup generated by $\ell$ and $2b$. 
The conductor of the full singularity $\Sigma_\infty$ is now given by
Gorenstein's formula \cite[Theorem 2]{gorenstein} and is

$$\sum_{1\le b \le m_\ell} \lbrace b(4m_\ell^2+4m_\ell-1)-2m_\ell-(2m_\ell+1)b^2\rbrace \cdot \sigma_{\infty,b}.$$

The full conductor ${\Cgot_\ell}$ is the sum of this plus the two
corresponding terms to the isomorphic   singularities $\Sigma_0$ and
$\Sigma_1$. 
The degree $\deg ({\Cgot_\ell})$ of ${\Cgot_\ell}$ is
$2m_\ell(2m_\ell^3+4m_\ell^2-2m_\ell-1)$. So we set 
$\delta =m_\ell(2m_\ell^3+4m_\ell^2-2m_\ell-1)$
and we check that $g_a=g_\ell+\delta$.

Now let $p\not \in \{2,3,\ell \}$ be a prime.  Let $\CC_p$ be the
(complete,  algebraically closed)  field of $p$-adics and $\Fpb$ its
residue
field. We embed $\Qb$ in  $\CC_p$ and also in $\CC$. In particular
$\zeta_\ell=\exp(\frac{2i\pi}{\ell})$
and $2^{\frac{1}{\ell}}$ are well defined as  $p$-adic numbers.
We observe that in the calculations above,  all  coefficients
belong to  $\ZZ [\frac{1}{6\ell}, \zeta_\ell, 2^{\frac{1}{\ell}}]$.
More precisely, the curves $C_\ell$ and $\cX_\ell$ are defined over $\ZZ [\frac{1}{6\ell}]$.
We write  $C_\ell \bmod p = C_\ell/\Fp = C_\ell\otimes_{ \ZZ
  [\frac{1}{6\ell}]}\Fp$ for the reduction of $C_\ell$ modulo $p$,  and
  define similarly $\cX_\ell\bmod p$.
We write similarly  
$\sigma_{\infty, b}\bmod   p$
and $\mu_{\infty, a}\bmod   p$.

We deduce the following.

\begin{lemma}[Computing $C_{\ell}$  and resolving its singularities]\label{lemma:computingCl}
There exists a deterministic algorithm that given a prime $\ell\ge 5$
and a prime  $p\not \in \{
2,3,\ell\}$ and a finite field  $\FF_q$ with characteristic  $p$ such that $\zeta_\ell \bmod p$ and
$2^{\frac{1}{\ell}}\bmod p$  belong to $\Fq$, computes
the equation $\cT_\ell(\lambda,x)$ modulo $p$ and   the expansions
of all $x_{a,b}$ as series in $\lambda^{\frac{-1}{\ell}}$ with
coefficients in $\Fq$, in time polynomial in $\ell$, $\log q$ and the
required  $\lambda^{\frac{-1}{\ell}}$-adic accuracy.  
\end{lemma}

\section{Another family  of modular curves}\label{section:modularcurves2}

In this section we consider another family of modular curves for which we
can
easily provide and study  a plane model. 
This family will be useful 
in the calculation of modular  representations as sketched in the next section.
Let $\ell > 5$ be a prime.  This time we set  $\cX_\ell=X_1(5\ell)$ the
moduli of elliptic curves with 
one  point of order $5\ell$. The
genus of $\cX_\ell$ is $g_\ell=\ell^2-4\ell+4$.
 We first describe a homogeneous singular plane  model $C_\ell$ for
this curve. We then  enumerate the geometric points on  $\cX_\ell$ above every 
singularity of $C_\ell$ and  provide series expansions for affine
coordinates at every such branch.
Finally, for $p\not\in\{2,3,5,\ell\}$ a prime integer,
we recall how to compute the zeta function of the function
field $\Fp(\cX_\ell)$. All this will be useful in section 
\ref{section:ramanujan} where we apply theorem \ref{theorem:mainjacobi}
to the curve $\cX_\ell$.

Let $b$ be an indeterminate and form the  elliptic curve $E_b$ in Tate
normal form 
with equation $y^2+(1-b)xy-by=x^3-bx^2$. The point $P=(0,0)$ has order
$5$ and its multiples are  $2P=(b,b^2)$, $3P=(b,0)$, $4P=(0,b)$.
The multiplication by $\ell$
isogeny
induces a degree $\ell^2$ rational function on $x$-coordinates:
$x\mapsto \frac{\cN(x)}{\cM(x)}$ where $\cN(x)$ is a monic degree $\ell^2$
polynomial in $\QQ(b)[x]$. Recursion formulae for division polynomial
(see \cite{enge} section 3.6) provide a quick algorithm for computing this
polynomial, and also show that the coefficients actually lie in
$\ZZ[b]$. If $\ell$ is congruent to $\pm 1$ modulo $5$ then $\ell P=\pm
P$ and $x$ divides $\cN(x)$. Otherwise $\cN(x)$ is divisible by $x-b$.

    Call $\cT_\ell(b,x)$ the
quotient of $\cN(x)$ by $x$ or $x-b$, accordingly.
 This  is 
 a monic  polynomial in $\ZZ[b][x]$ with 
degree
$\ell^2-1$ in $x$. As a polynomial in $x$ we have

$$\cT_\ell(b,x)=\sum_{0\le k \le
 \ell^2-1}a_{\ell^2-1-k}(b)x^k$$
\noindent where  $a_0(\lambda)=1$.  We call
 $d$ be the total degree
of $\cT_\ell$.

As in the previous section, we check that $\cT_\ell$
is irreducible in $\bar \QQ (b)[x]$ and $\QQ(\cX_\ell)$
is the splitting field of $\cT_\ell$ over
$\QQ(b)$.
Let  $C_\ell\subset \PP^2$ be the  projective curve
made of the zeroes of the homogeneous 
polynomial $\cT_\ell(\frac{B}{Y},\frac{X}{Y})Y^{d}$.

We set 
$$j=j(b)=\frac{(b^4-12b^3+14b^2+12b+1)^3}{b^5(b^2-11b-1)}.$$

Let $\sqrt{5}\in \CC$ be the positive square root of $5$ and 
let $\zeta_5=\exp(\frac{2i\pi}{5})$. Let   $s=\frac{11+5\sqrt{5}}{2}$ and $\bar s$ be the two roots of
$b^2-11b-1$.
The forgetful map $X_1(5\ell)\rightarrow X_1(5)$ is unramified
except at $b\in \{0,\infty,s,\bar s\}$.
For every point $P$ on $\cX_\ell$ such that $b(P) \not \in
\{0,s,\bar s,\infty \}$,  the function $b -b (P)$ is a
uniformizing parameter at $P$. 

Let $\cU$ be the affine open set with equation $YB(B^2-11BY+Y^2)\not =
0$. Every point  on  $C_\ell\cap \cU$  is smooth and all points on
$\cX_\ell$ above  points in
$C_\ell - \cU$ are cusps in the modular sense (i.e. the modular
invariant at these points is infinite).

In order to desingularize $C_\ell$ at a given  cusp, we shall construct an
isomorphism between the Tate $q$-curve and the completion of $E_b$
at this  cusp.
We call $A_\infty$, $A_0$, $A_s$, $A_{\bar s}$ the points on $X_1(5)$
corresponding to the values $\infty$, $0$, $s$ and $\bar s$ of $b$.
We first study the situation locally at $A_\infty$.
A local parameter is $b^{-1}$ and $j^{-1}=b^{-5}+25b^{-6}+\cdots$.

We introduce Tate's $q$-parameter,
defined implicitly by

$$j=\frac{1}{q}+744+196884q+\cdots$$
\noindent  so

\begin{eqnarray*}
q&=&j^{-1}+744j^{-2}+750420j^{-3}+\cdots\\
&=& b^{-5}+25b^{-6}+\cdots
\end{eqnarray*}
\noindent  and
we fix an embedding of the local field at 
$A_\infty$ inside the field of Puiseux series
$\CC\{\{q\}\}$ by setting $b^{-1}=q^\frac{1}{5} - 5q^\frac{2}{5}+\cdots$.

We set $x'=36x+3(b^2-6b+1)$ and $y'=108(2y+(1-b)x-b)$ and  find the reduced
Weierstrass equation 
$$y'^2=x'^3-27(b^4-12b^3+14b^2+12b+1)x' +54 (b^2+1  ) ( b^4 -18b^3 +74 b^2+18b+1  ).$$

We want to compare the latter curve and the Tate  curve with equation
$$y''^2=x''^3-\frac{E_4(q)}{48}x''+\frac{E_6(q)}{864}$$ where
$E_4(q)=1+240q+\cdots$ and $E_6(q)=1-504q+\cdots$. See
\cite[Theorem 10.1.6]{husemoller}.

From the classical (see \cite[Proposition 7.1]{schoof}) identities 

$$\left( \frac{qdj}{dq}  \right)^2=j(j-1728)E_4$$
\noindent 

$$\left(
\frac{qdj}{dq}  \right)^3=-j^2(j-1728)E_6$$
\noindent 
 we deduce  
$$\left(
\frac{qdb}{dq}  \right)^2=
\frac{b^2(b^2-11b-1)^2E_4}{25(b^4-12b^3+14b^2+12b+1)}$$
\noindent  and 
$$\left(
\frac{qdb}{dq}  \right)^3=-\frac{b^3(b^2-11b-1)^3E_6}{125(b^2+1)(b^4-18b^3+74b^2+18b+1)}.$$

We deduce the isomorphism $x'=\gamma ^2x''$  and $y'=\gamma^3y''$ with

$$\gamma^2=-\frac{36b(b^2-11b-1)dq}{5qdb}.$$

The point $P$ has $(x,y)$ coordinates equal to $(0,0)$. So

$$x''(P)=3(b^2-6b+1)/\gamma^2=\frac{1}{12}+b^{-2}+11b^{-3}+\cdots= \frac{1}{12}+q^{\frac{2}{5}}+O(q^{\frac{3}{5}}).$$

Since on  the Tate  curve we have

\begin{equation}\label{equation:expa}
x''(w, q)=\frac{1}{12}+\sum_{n\in \ZZ}\frac{wq^n}{(1-wq^n)^2}
-2\sum_{n\ge 1}\frac{nq^n}{1-q^n}
\end{equation}
\noindent we deduce that  $w(P)=q^{\pm \frac{2}{5}} \bmod <q>$. We may
take either  sign in the exponent because we may choose  any of  the two
isomorphisms corresponding to either possible values for $\gamma$. We
decide that $w(P)=q^{\frac{2}{5}} \bmod <q>$.
Set  $\zeta_\ell=\exp(\frac{2i\pi}{\ell})$.
For $\alpha$ and $\beta$ integers such that $0\le \alpha,\beta\le \ell-1$    we set
$w=\zeta_\ell^\alpha q^{\frac{\beta}{\ell}}q^{\frac{2}{5\ell}}$ in the  expansion (\ref{equation:expa})
and find

$$x''_{\alpha,\beta}=\frac{1}{12}+\zeta_\ell^\alpha q^{\frac{\beta}{\ell}}q^{\frac{2}{5\ell}}(1+O(q^{\frac{1}{5\ell}}))$$
\noindent if $0\le \beta \le \frac{\ell -1}{2}$ and
$$x''_{\alpha,\beta}=\frac{1}{12}+\zeta_\ell^{-\alpha}q^{\frac{\ell-\beta}{\ell}-\frac{2}{5\ell}}(1+O(q^{\frac{1}{5\ell}}))$$
\noindent if $\frac{\ell +1}{2}\le \beta \le \ell-1$.

Since

$$x_{\alpha,\beta}=(\gamma^2
x''_{\alpha,\beta}-3(b^2-6b+1))/36$$
\noindent and  $\gamma^2= 36b^{2} - 216b -
396+O(b^{-1}) = 36q^\frac{-2}{5} +144q^\frac{-1}{5} +144+\cdots$ we
deduce that 

$$x_{\alpha,\beta}+1=\zeta_\ell^\alpha q^{\frac{\beta}{\ell}+\frac{2}{5\ell}-\frac{2}{5}}(1+O(q^{\frac{1}{5\ell}}))$$
\noindent if $0\le \beta \le \frac{\ell -1}{2}$  and

$$x_{\alpha,\beta}+1=\zeta_\ell^{-\alpha} q^{\frac{\ell-\beta}{\ell}-\frac{2}{5\ell}-\frac{2}{5}}(1+O(q^{\frac{1}{5\ell}}))$$
\noindent if $\frac{\ell +1}{2}\le \beta \le \ell-1$.

 In particular, the degree of $\cT_\ell(b,x)$ in $b$ is $\le 2(\ell^2-1)$.

For $0\le \alpha <\ell$ and $0\le \beta <\ell$ we set $\tilde \alpha =
5\alpha \bmod \ell$ and $\tilde \beta =5\beta+2\bmod \ell$.
If  $\tilde \beta$ is non-zero,  the local monodromy
 group permutes cyclically
the $\ell$ roots  $x_{\alpha, \beta}$ for $0\le \alpha < \ell$. We
call
$\sigma_{\infty,\tilde \beta}$ the corresponding branch on $\cX_\ell$.
 On
the other hand, if $\beta = \frac{-2}{5}\bmod \ell$ then $\tilde \beta
=0\bmod \ell$ and every $x_{\alpha,\frac{-2}{5}\bmod \ell}$ is fixed
by the local monodromy group. 
We observe that $x_{0,\frac{-2}{5}\bmod \ell}$ is either $b$ or $0$
and is not a root of $\cT_\ell(b,x)$.   For
$\tilde \alpha $ a non-zero residue modulo $\ell$, we denote  by $\mu_{\infty,
  \tilde \alpha}$ the
 branch on $\cX_\ell$ corresponding to $x_{\alpha,\frac{-2}{5}\bmod \ell}$.

So we have
$\ell-1$
unramified points  on $\cX_\ell$ above $A_\infty$ and $\ell -1$
ramified points
with
ramification index $\ell$.

The coefficients in  all the series expansions above are in
 $ \ZZ [\frac{1}{30},\zeta_\ell]$. The coefficients of $\cT_\ell(b, x)$
 are
in  $\ZZ$. From the discussion above we deduce the following.


\begin{lemma}[Computing $C_{\ell}$  and resolving its singularities, I]\label{lemma:computingCl2}
There exists a deterministic algorithm that given a prime $\ell\ge 7$
and a prime  $p\not \in \{
2,3,5,\ell\}$ and a finite field  $\FF_q$ with characteristic  $p$ such that $\zeta_\ell \bmod p$ 
  belongs to $\Fq$, computes
the equation $\cT_\ell(b,x)$ modulo $p$ and   the expansions
of all $x_{\alpha,\beta}$ as series in $b^{-\frac{1}{\ell}}$ with
coefficients in $\Fq$, in time polynomial in $\ell$, $\log q$ and the
required  $b^{\frac{-1}{\ell}}$-adic accuracy.  
\end{lemma}

In appendix \ref{app:A} we give
a  few lines of GP-PARI code (see \cite{pari}) that compute
these expansions.

\medskip 
We now study
the  singular points above $A_0$. A local parameter at
$A_0$
is $b$ and $j^{-1}=-b^5+25b^6+\dots$ so $q=-b^5+25b^6+\dots$ and we
fix an embedding of the local field at $A_0$ inside $\CC\{\{q\}\}$ by setting
$b=-q^{\frac{1}{5}}+5q^{\frac{2}{5}}+\dots$. From
$\gamma^2=36-216q^{\frac{1}{5}}+\dots$
we deduce that the coordinate $x''(P)$ of the $5$-torsion point $P$ is
$x''(P)=\frac{1}{12}+q^{\frac{1}{5}}+O(q^{\frac{2}{5}})$ so the
parameter $w$ at $P$ can be taken to be  $w(P)=q^{\frac{1}{5}}\bmod <q>$ this time.
For $\alpha$ and $\beta$ integers such that $0\le \alpha,\beta\le \ell-1$    we set
$w=\zeta_\ell^\alpha q^{\frac{\beta}{\ell}}q^{\frac{1}{5\ell}}$ in the  expansion (\ref{equation:expa})
and we finish as above.

Now, a local parameter at
$A_s$
is $b-s$ and $j^{-1}=(\frac{1}{2}-\frac{11\sqrt{5}}{50})(b-s)+O((b-s)^2)$ so $q=(\frac{1}{2}-\frac{11\sqrt{5}}{50})(b-s)+O((b-s)^2)$ and we
fix an embedding of the local field at $A_s$ inside $\CC\{\{q\}\}$ by setting
$b-s=\frac{125+55\sqrt{5}}{2}q+O(q^2)$. 
We deduce that the coordinate $x''(P)$ of the $5$-torsion point $P$ is
$x''(P)=\frac{1}{12}+\frac{w}{(1-w)^2}+O(q)$ where $w=\exp(\frac{4i\pi}{5})=\zeta_5^2$
so the
parameter $w$ at $P$ can be taken to be  $w(P)= \zeta_5^{2}\bmod <q>$ this time.

Altogether we have proved the following.

\begin{lemma}[Computing $C_{\ell}$  and resolving its singularities, II]\label{lemma:computingCl22}
There exists a deterministic algorithm that given a prime $\ell\ge 7$
and a prime  $p\not \in \{
2,3,5,\ell\}$ and a finite field  $\FF_q$ with characteristic  $p$
such that $\zeta_\ell \bmod p$ 
and $\zeta_5 \bmod p$
  belong to $\Fq$, computes
the equation $\cT_\ell(b,x)$ modulo $p$ and    expansions (with 
coefficients in $\Fq$)
at every singular branch of $C_\ell$  in time polynomial in $\ell$, $\log q$ and the
required  number of significant terms in the expansions.  
\end{lemma}

In order to apply theorem
 \ref{theorem:mainjacobi}
to the curve $\cX_\ell$, we shall
also need the following result due to Manin, Shokurov, Merel and
Cremona \cite{manin, merel, cremona, frey}.  

\begin{lemma}[Manin, Shokurov, Merel, Cremona]\label{lemma:manin2}
For $\ell$ a prime and  $p\not\in \{5,\ell\}$ another prime, the
zeta function of $\cX_\ell \pmod p$ can be  computed in deterministic 
polynomial time
in $\ell$ and  $p$.
\end{lemma}

We first  compute the action of the Hecke
operator $T_p$ on the space of Manin symbols for
the congruence group $\Gamma_1(5\ell)$ associated with  $\cX_\ell$.
Then, from the Eichler-Shimura identity $T_p=F_p+p<p>/F_p$ we deduce the characteristic polynomial
of the Frobenius $F_p$. \hfill $\Box$

In appendix \ref{app:B} we give
a  few lines of Magma code (see \cite{magma}) that compute
the zeta function of $X_1(5\ell)/\Fp$.

\section{Computing the Ramanujan subspace 
over $\Fp$}\label{section:ramanujan}

This section explains the connection  between the methods given here
and  Edixhoven's program for computing coefficients of modular forms.
Recall the definition of the Ramanujan arithmetic $\tau$ function,
related to the sum  expansion  of the discriminant form:

$$\Delta(q)=q\prod_{k\ge 1}(1-q^k)^{24}=\sum_{k\ge 1}\tau(k)q^k.$$

We call $\TT\subset \End(J_1(\ell)/\QQ)$ the algebra  of endomorphisms of $J_1(\ell)$
generated by the  Hecke  operators $T_n$ for all  integers $n\ge 2$. 
Following Edixhoven \cite[Definition 10.9]{arxiv}  we state
the 

\begin{definition}[The Ramanujan ideal]\label{definition:vl}
Assume $\ell \ge 13$ is a prime.   We denote by  $\mgot$ the maximal
ideal
in  $\TT$ generated by $\ell$ and  the $T_n-\tau(n)$.
The   subspace $J_1(\ell)[\mgot]$ of the $\ell$-torsion of $J_1(\ell)$
 cut out by all $T_n-\tau(n)$ is called the Ramanujan subspace at $\ell$
and denoted   $V_\ell$.
\end{definition}

This $V_\ell$ is a $2$-dimensional vector space over $\Fl$ and
for $p\not = \ell$
the characteristic polynomial of the Frobenius endomorphism
$F_p$ on it is $X^2-\tau(p)X+p^{11}\bmod \ell$.

In this section, we address the problem of computing 
$\mgot$-torsion divisors on modular curves over some  extension field
$\Fq$ of $\Fp$ for $p\not = \ell$. The definition field $\Fq$ for such divisors  can be predicted from the
characteristic polynomial of $F_p$ on $V_\ell$. So the strategy is to
pick random $\Fq$-points in the
$\ell$-torsion of the  jacobian
$J_1(\ell)$ and to project them onto $V_\ell$ using Hecke operators.

In section \ref{section:modularcurves2} we have defined the modular
curve $\cX_\ell=X_1(5\ell)$
and the degree $24$ covering $\phi : \cX_\ell
\rightarrow X_1(\ell)$ of $X_1(\ell)$. We prefer $\cX_\ell$  to $X_1(\ell)$
because we are able to construct  a natural and convenient 
plane model  for it. 
The  covering map $\phi :  \cX_\ell\rightarrow X_1(\ell)$
corresponds to forgetting the $5$-torsion structure. It induces
two morphisms $\phi ^* : J_1(\ell)\rightarrow \cJ_\ell$ and $\phi_* :
\cJ_\ell\rightarrow J_1(\ell)$ such that the composite map 
$\phi_*\circ \phi^*$ is   multiplication by $24$ in  $J_1(\ell)$. We write
$\phi_*\circ \phi^*=[24]$. 
Thus the curve $\cX_\ell$ provides a convenient computational
model for the
group of $\Fq$-points of  the jacobian of  $X_1(\ell)$.


We denote by $\cA_\ell \subset \cJ_\ell$ the image of 
$\nu = \phi^*\circ\phi_*$. This is a subvariety of $\cJ_\ell$ isogenous to
$J_1(\ell)$. 
The restriction of $\nu$ to $\cA_\ell$ is multiplication by $24$.
The maps $\phi^*$  and $\phi_*$ induce Galois equivariant
bijections between
the
$N$-torsion subgroups $J_1(\ell)[N]$ and $\cA_\ell[N]$ for every integer $N$
which is prime to
$6$.

We call $W_\ell\subset \cA_\ell\subset \cJ_\ell$ the image of the Ramanujan subspace 
by $\phi^*$. 
We  choose an integer $k$ such that $24k$ is congruent
to $1$ modulo $\ell$, and  set  $\hT_n=[k]\circ \phi^*\circ T_n\circ
\phi_*$, for every $n$. We  notice that $\hT_n\circ \phi^*=\phi^*\circ T_n$
on $J_1(\ell)[\ell]$. This
way, the map 
$\phi^* : J_1(\ell)\rightarrow
\cJ_\ell$ induces   a Galois equivariant bijection of Hecke modules
between   $J_1(\ell)[\ell]$ and  $\cA_\ell[\ell]$, and $W_\ell=\phi^*(V_\ell)$
 is the subspace in $\cA_\ell[\ell]$  cut out by all
$\hT_n-\tau(n)$. So $W_\ell$ will also be called the Ramanujan
subspace at $\ell$
whenever there is no risk of confusion.
We notice that $\phi^*$, $\phi_*$, $T_n$, and $\hT_n$ can be seen as
correspondences
as well as morphisms between  jacobians, and we state   the  following.

\begin{lemma}[Computing the Hecke action]\label{lemma:hecke}
Let  $\ell$ and  $p$ be primes such that $p\not \in \{2,3,5,\ell\}$.
 Let $n\ge 2$
be an integer.  Let $q$ be a
power of $p$ and let $D$ be an effective $\Fq$-divisor of degree
$\deg (D)$ on $\cX_\ell \pmod p$. The divisors $\phi^*\circ \phi_*(D)$ and $\phi^*\circ T_n\circ \phi_*(D)$
  can be computed in  polynomial time in
  $\ell$, $\deg(D)$, $n$ and $\log q$.
\end{lemma}

If $n$ is prime to $\ell$, we define
the Hecke operator $T(n,n)$ as an element in the ring of correspondences 
on $X_1(\ell)$ tensored by $\QQ$. See \cite[VII, \S 2 ]{langmf}.
From \cite[VII, \S 2, Theorem 2.1]{langmf} we have 
$T_{\ell^i}=(T_\ell)^i$ and $T_{n^{i}}=T_{n^{i-1}}T_n-nT_{n^{i-2}}T(n,n)$
if $n$ is prime and $n\not =\ell$.
And of course $T_{n_1}T_{n_2}=T_{n_1n_2}$ if $n_1$ and $n_2$ are
coprime.
So it suffices to explain how to compute $T_\ell$ and also 
$T_n$ and $T(n,n)$ for $n$ prime and
$n \not =\ell$.

Let  $x=(E,u)$ be  a point on $Y_1(\ell)\subset X_1(\ell)$
 representing an elliptic
curve $E$ with one $\ell$-torsion point $u$. 
Let $n$ be an integer. The Hecke operator $T_n$  maps
$x$ onto the sum of all $(E_I,I(u))$, where $I : E\rightarrow
E_I$ runs over the set of all  isogenies of degree $n$ from
$E$ such that $I(u)$ still has order $\ell$.  
If $n$ is prime to
$\ell$, the Hecke operator $T(n,n)$  maps $x$ onto
$\frac{1}{n^2}$ times $(E,nu)$. 
So we can compute the action of these  Hecke correspondences on points
$x=(E,u)$ using V{\'e}lu's formulae \cite{velu}.

There remains to treat the case of cusps.
We call $\sigma_{\tilde \beta}$ for $1\le \tilde \beta\le
\frac{\ell-1}{2}$ and $\mu_{\tilde \alpha}$ for
$1\le \tilde{\alpha}\le \frac{\ell -1}{2}$ the cusps on $X_1(\ell)$ images by
$\phi$ of the 
$\sigma_{\infty,\tilde \beta}$ and $\mu_{\infty,\tilde \alpha}$.
To every cusp one can associate a  set of Tate curves
with $\ell$-torsion  point (one Tate curve  for every branch
at this cusp).

For example the Tate curves at  $\sigma_\tbeta$
are the Tate curves $\CC^*/q$ with $\ell$-torsion
point $w=\zeta_\ell^\star q^{\frac{ \tbeta}{\ell}}$ where the star runs
over the set of all
 residues modulo $\ell$. There are $\ell$ branches at each such cusp.

Similarly, the Tate curves at $\mu_\talpha$
are the Tate curves $\CC^*/q$ with $\ell$-torsion
point $w=\zeta_\ell^{ \talpha}$. One single branch here:
no ramification.

For $n$ prime and  $n\not = \ell$ we have

$$T_n(\sigma_{\tilde\beta})=\sigma_{\tilde
 \beta}+n\sigma_{n\tilde\beta}$$
\noindent  and
$$T_n(\mu_{\tilde \alpha})=n\mu_{\tilde\alpha}+\mu_{n\tilde\alpha},$$
\noindent 
where $n\tilde \alpha$ in $\mu_{n\tilde \alpha}$  (resp. $n\tilde\beta$ in $\sigma_{n\tilde\beta}$)
 should be understood as a class in
$(\ZZ/\ell\ZZ)^*/\{ 1,-1\}$.

Similarly
$$T_\ell(\sigma_{\tilde\beta})=\sigma_{\tilde\beta}+2\ell\sum_{1\le
 \tilde\alpha\le \frac{\ell -1}{2}}\mu_{\tilde\alpha}$$
\noindent  and
$$T_\ell(\mu_{\tilde\alpha})=\ell\mu_{\tilde\alpha}.$$

And of course, if $n$ is prime to $\ell$, then   $T(n,n)(\sigma_\tbeta)=
\frac{1}{n^2}\sigma_{n\tbeta}$ and $T(n,n)(\mu_\talpha)=
\frac{1}{n^2}\mu_{n\talpha}$.

All together, one can compute the effect of $T_n$ on
cusps for all $n$. 
For the sake of completeness, we also give the action of the diamond
operator $<\!\! n\!\!>$ on cusps. If $n$ is prime to $\ell$ then
$<\!\!n\!\!>(\sigma_\tbeta)=\sigma_{n\tbeta}$ and $<\!\!n\!\!>(\mu_\talpha)=\mu_{n\talpha}$.

\hfill $\Box$

We can now state the  following.

\begin{theorem}\label{theorem:computingltorsionmodp}
There is a probabilistic (Las Vegas)  algorithm that on input  a prime
$\ell\ge 13$ 
and a prime $p\ge 7$ such that $\ell\not = p$,
computes
the Ramanujan subspace $W_\ell=\phi^*(V_\ell)$  inside the $\ell$-torsion of
the jacobian of $\cX_\ell / \Fp$. The answer is given as a list
of $\ell^2$ degree $g_\ell$ effective divisors on $\cX_\ell$, the
first one being the origin $\omega$.
The algorithm runs in probabilistic polynomial time in $p$
and $\ell$.
\end{theorem}
Lemma \ref{lemma:computingCl22} gives us a  plane model for $\cX_\ell\pmod p$ and a resolution of its
singularities. From lemma \ref{lemma:manin2} we obtain the zeta function of
$\cX_\ell \pmod p$.  
The characteristic polynomial of $F_p$ on the Ramanujan space $V_\ell$
is $X^2-\tau(p)X+p^{11}\bmod \ell$.
So we compute $\tau (p) \pmod \ell$ using the expansion of the
discriminant
form. We deduce some small enough
 field of decomposition $\Fq$ for $V_\ell \pmod p$.
We then apply  theorem \ref{theorem:mainjacobi} and obtain
a basis for the $\ell$-torsion in the Picard group of $\cX_\ell
/\Fq$. The same theorem allows us to compute the matrix of the
endomorphism $\nu = \phi^*\circ\phi_*$
in this basis. We deduce  a basis for the image $\cA[\ell](\Fq)$ of $\nu$.
Using  theorem \ref{theorem:mainjacobi} again, we now write down the matrices of the Hecke operators $\hT_n$ in this
basis
for all $n < \ell^2$. It is then a matter of linear algebra to compute a
basis
for the intersection of the kernels of all $\hT_n -\tau (n)$ in
$\cA[\ell](\Fq)$. The algorithm is Las Vegas rather than Monte-Carlo
because we 
can check the result,  the group  $W_\ell$ having known cardinality $\ell^2$.
\hfill $\Box$

\begin{remark}
In the above theorem, one may impose an origin $\omega$ rather than
letting
the algorithm choose it. For example, following work
by Edixhoven in \cite[Section 12]{arxiv}, one may choose  as origin a well designed linear
combination
of the cusps. Such  an adapted choice of the origin may ensure that
the $\ell^2-1$ divisors representing the non-zero classes in $W_\ell$
are unique in characteristic zero and thus remain unique modulo $p$
for all but finitely many primes $p$.
\end{remark}

\section{The semisimple non-scalar case}\label{section:semi}

In this section we present a simplified algorithm
for computing the Ramanujan subspace $V_\ell$ modulo $p$,  that
applies when the Frobenius action on it is semisimple
and non-scalar or equivalently when  $\tau(p)^2-4p^{11}$ is not divisible by $\ell$.
The main idea is to associate a divisible group with $V_\ell$.

For every integer $n\ge 2$ we call $A_n(X)\in \ZZ[X]$ the
characteristic 
polynomial of $T_n$ acting on weight
$2$ modular forms for $\Gamma_1(\ell)$. We factor 

$$A_n(X)=B_n(X)(X-\tau(n))^{e_n}$$
\noindent 
in $\Fl[X]$ 
with
$B_n(X)$ monic and $B_n(\tau(n))\not = 0\in \Fl$. 
For every integer $k\ge 1$ this  polynomial factorization   lifts
modulo $\ell^k$ as 

$$A_n(X)= B_{n,k}(X) C_{n,k}(X)\pmod {\ell^k}.$$

We call $\Pi_{k} : J_1(\ell)[\ell^k]\rightarrow J_1(\ell)[\ell^k]$
the composite map of  all 
$B_{n,k}(T_n)$ for all integers $n$  such that $2\le
n < \ell^2$. We observe that $\Pi_{{k+1}}$ coincides with
$\Pi_{k}$ on  $J_1(\ell)[\ell^k]$.  So we have defined a map 
$\Pi : J_1(\ell)[\ell^\infty]\rightarrow J_1(\ell)[\ell^\infty]$.

We have the following.

\begin{lemma}[The  Ramanujan modules]\label{lemma:G}
For $k\ge 1$ an integer, we denote by $\GG_k$ the
subgroup of  $J_1(\ell)[\ell^k]$ consisting of elements killed
by some power of $\mgot$.  Let $\GG$ be the union of all $\GG_k$.
The group 
$\GG_k$ is the image $\Pi_k(J_1(\ell)[\ell^k])$ of the
$\ell^k$-torsion by $\Pi_k$. It is killed by $\mgot^{2kg(X_1(\ell))} $ and the
restriction of $\Pi_{k}$ to $\GG_k$ is a bijection. Further 
$\GG_{k+1}[\ell^k]= \GG_k = \ell \GG_{k+1}$. 
The 
$(\ZZ/\ell^k\ZZ)$-module
$\GG_k$ is free. We call it the Ramanujan module.
\end{lemma}

We show that for every integer $n\ge 2$,
the  restriction of $B_{n,k}(T_n)$ to $\GG_k$
 is a bijection. It suffices to show
injectivity. Assume $B_{n,k}(T_n)$ restricted to $\GG_k$ is not
 injective. There is a non-zero $\ell$-torsion  element  $P$  in its
 kernel.
This $P$ is killed by $(T_n-\tau(n))^m \pmod \ell$ for some integer $m$. It
is also killed by $B_{n}(T_n)\pmod \ell$. Since these two polynomials
are coprime, $P$ is zero, contradiction.

So $\Pi_k$ is an automorphism of $\GG_k$. In particular $\GG_k\subset
\Pi_k (J_1(\ell)[\ell^k])$.
We set $\II_k=\Pi_k (J_1(\ell)[\ell^k])$ and we prove the converse
inclusion $\II_k\subset \GG_k$.
For every integer $n$ between $2$ and $\ell^2$,
the restriction of $T_n$ to $\II_1$ is killed by $(X-\tau(n))^{e_n}$.
Since the Hecke algebra is generated by these $T_n$ and is
commutative, its image in $\End(\II_1)$ is triangulisable
\footnote{If $K$ is a field and $V$ a $K$-vector space, we write
$\cL(V)$ for the algebra of linear maps from $V$
to itself. Let $A$ be a subset of $\cL(V)$.  We say that $A$ is 
triangulisable if  there exists a
  basis $\cB$ of $V$
 such that  the matrix of
 every element 
  in $A$  with respect to $\cB$ is upper triangular.} and consists
of matrices with a single eigenvalue. We
 deduce that for every integer $n$ the restriction of $T_n$
to $\II_1$ has a single eigenvalue (namely $\tau(n)\pmod
\ell$). 
Because the dimension  of $\II_1$ as a $\Fl$-vector space is
 $\le 2g(X_1(\ell))$ we deduce that $\II_1$ is killed by 
$\mgot^{2g(X_1(\ell))}$. So $\II_1 = \GG_1$ is killed by
$\mgot^{2g(X_1(\ell))}$.

For every  integer $n$ between $2$ and $\ell^2$,
the restriction of $T_n$ to $\II_k[\ell]$ is killed by $C_{n,k}(X)$
which is congruent to $(X-\tau(n))^{e_n}$ modulo $\ell$. 
So 
$\II_k[\ell]$ is killed by $(T_n-\tau(n))^{e_n}$ and by 
$\mgot^{2g(X_1(\ell))}$. So any
morphism in $\mgot^{2kg(X_1(\ell))}$  kills
$\II_k[\ell^k]=\II_k$. So  $\II_k$ is
killed
by  $\mgot^{2kg(X_1(\ell))}$ and $\II_k=\GG_k$.

It is clear that $\ell\GG_{k+1}\subset \GG_k$. Conversely if
$P=\Pi_k(Q)$ and  $Q$ is  $\ell^k$-torsion then let $R$ such that
$\ell R=Q$ and $S=\Pi_{k+1}(R)$. Then $S$ is in $\II_{k+1}=\GG_{k+1}$ and
$\ell S= \Pi _{k+1}(Q)= \Pi_k(Q)=P$. So $\ell \GG_{k+1}=\GG_k$.
From $\GG_{k+1}[\ell^k]=\ell\GG_{k+1}$ we deduce that $\GG_{k+1 }$ is
a free $(\ZZ/\ell^{k+1}\ZZ)$-module.  \hfill $\Box$

\medskip
We now study the Galois action on this divisible group.
Let $p\not =\ell$ be a prime. We regard  $J_1(\ell)$ as a variety over
the finite field $\Fp$. The Ramanujan module
$\GG=J_1(\ell)[\mgot^\infty]$ is then an
$\ell$-divisible group inside $J_1(\ell)[\ell^\infty]$ in the sense
of definition \ref{definition:divisible}.
According to the Eichler-Shimura identity
$F_p^2-T_pF_p+p<p>=0$.
The diamond operator $<\!\!p\!\!>\in \TT$ has a
unique eigenvalue on $\GG_1$, namely  $p^{10}\pmod \ell$. Since
$F_p$ commutes with $\TT$, the algebra generated by
$\TT$ and $F_p$ is triangulisable\footnotemark[1]
in $\GL(\GG_1\otimes_\Fl\bar\FF_\ell)$. So any eigenvalue of $F_p$ on $\GG_1$ is
killed
by $X^2-\tau(p)X+p^{11}\pmod \ell $.
Let $\eta$ be an integer that kills the roots of 
 the  polynomial $X^2-\tau(p)X+p^{11} \pmod \ell$ in
 $\bar\FF_\ell^*$. For example one may take $\eta = \ell^2-1$.
As an endomorphism of $\GG_1$ one has
 $F_p^\eta = \Id +n$ where $n$ is nilpotent. Since the dimension
 of $\GG_1$ is $\le 2g(X_1(\ell))\le \ell^2$ one has $n^{\ell^2}=0$
 and $F_p^{\eta\ell^2}=\Id$.   So  $\GG_1$ splits completely over
 $\FF_{p^{\ell^2(\ell^2-1)}}$. 
As a consequence, $\GG_k$ splits completely over the extension of degree 
${(\ell^2-1)}\ell^{k+1}$ of $\FF_{p}$. 

\begin{lemma}[Galois action on the Ramanujan module]\label{lemma:G2}
If $p\not =\ell$ is a prime, then the Ramanujan module
$\GG=J_1(\ell)[\mgot^\infty]$  is a
divisible group inside $J_1(\ell)[\ell^\infty](\Fpb)$. Let $\eta$ be an
integer
that kills the roots of $X^2-\tau(p)X+p^{11}$ in $\bar
\FF_\ell^*$. For example  $\eta = \ell^2-1$. The $\ell^k$-torsion
$\GG_k=\GG[\ell^k]$  inside $\GG$ splits completely
 over the extension of degree $\eta
 \ell^{k+1}$ of $\FF_{p}$. 
\end{lemma}

For computational convenience we may prefer $\cX_\ell=X_1(5\ell)$ to
$X_1(\ell)$.  If this is the case, we embed $\GG$ inside the jacobian
$\cJ_\ell$ of $\cX_\ell$ using the map $\phi^*$. For the sake of
simplicity
we present the  calculations below  in the context of  $J_1(\ell)$
although they take place inside $\cJ_\ell$.

\medskip 

The knowledge of a non-zero element in $\GG_k$  sometimes
 suffices to
construct a basis of $V_\ell(\Fpb)$:

\begin{lemma}[The inert case]\label{lemma:projnonsplit}
Assume $X^2-\tau(p)X+p^{11}\pmod \ell$ is irreducible. 
Let $k\ge 1$ be an integer and $q=p^d$ a power of $p$. Given a non
zero element in $\GG_k(\FF_{q})$, one can compute a basis of
$V_\ell(\Fpb)$ in  polynomial time in $\log q$, $\ell$ and $k$.
\end{lemma}

Indeed, let $P\in \GG_k(\Fq)$ be non-zero. We replace $P$ by $\ell P$
until
we find a non-zero element in $\GG_1(\Fq)$. 
Given such a $P$ we can test whether it belongs
to $V_\ell$ by computing   $(T_n-\tau(n))x$  for all $2\le n\le
\ell^2$. If we only obtain zeroes this shows $P$ is in
$V_\ell$. Otherwise we replace $P$ by some non-zero $(T_n-\tau(n))P$
and test again. This process stops after $2g(X_1(\ell))$ steps at
most, and produces a non-zero element $P$ in $V_\ell(\Fq)$. Since $F_p$ has
no eigenvector in $V_\ell(\Fpb)$, the couple $(P,F_p(P))$ is a basis of
$V_\ell(\Fpb)$.\hfill $\Box$

\medskip 

So  assuming that  $\tau(p)^2-4p^{11}$ is not a square modulo $\ell$,
we have a simpler method to construct a basis for 
the Ramanujan module $V_\ell$ modulo $p$:

We set  $q=p^{(\ell^2-1)\ell^3}$. We have 
$\GG(\Fq)\supset \GG_2=\GG_2(\Fq)$. Set 
$N_q=\#J_1(\ell)(\Fq)=M_qL_q$ where $M_q$ is prime to $\ell$. 
This $N_q$ can be computed using Manin symbols as in lemma
\ref{lemma:manin2}.
Let $L_q=\ell^w$.
The image of $J_1(\ell)(\Fq)$ by the morphism $\psi = \Pi_w\circ
[M_q]$ contains $\GG_2(\Fq)$ and is in fact equal to  $\GG(\Fq)$. 
We check $\#\GG(\Fq)\ge \#\GG_2 \ge
\ell^4$. So at least one  of the elements in $J_1(\ell)(\Fq)$ given by  lemma \ref{lemma:generators}
has
a non-zero  image by $\psi$ for $\ell$  large enough.
We apply lemma \ref{lemma:projnonsplit} to this element and find a basis
for the Ramanujan module at $\ell$.

\medskip

We now assume the  polynomial $X^2-\tau(p)X+p^{11}\bmod \ell$ has 
two distinct roots
$a\bmod \ell$ and $b\bmod \ell$. So
$(F_p-a)^{2g(X_1(\ell))}(F_p-b)^{2g(X_1(\ell))}$ kills $\GG_1$.
 Since
$\GG_1=\GG_k[\ell]$ we deduce
that $(F_p-a)^{2kg(X_1(\ell))}(F_p-b)^{2kg(X_1(\ell))}$ kills $\GG_k$. 

This leads us to the following definition.

\begin{definition}[Split Ramanujan modules]
Assume $X^2-\tau(p)X+p^{11}\bmod \ell$ has two distinct 
roots $a\bmod \ell$ and $b\bmod \ell$ where $a$ and  $b$  are
integers. Let $\mgot_a$ be the ideal in $\TT[F_p]$ generated 
by $\ell$, all $T_n-\tau(n)$ and $F_p-a$. Let
$V_{\ell,a}=J_1(\ell)[\mgot_a]\subset V_\ell$
be the eigenspace associated with $a$.
For $k\ge 1$ an integer, we denote by $\GG_{k,a}$ the
subgroup of  $J_1(\ell)[\ell^k]$ consisting of elements killed
by some power of $\mgot_a$. 
Let $\Pi_{k,a}$ the composition of $\Pi_k$ and
$(F_p-b)^{2kg(X_1(\ell))}$.
We denote by $\GG_a$ the union of all $\GG_{k,a}$.
\end{definition}

We have the following.

\begin{lemma}[Properties of split Ramanujan modules]\label{lemma:Ga}
For every integer $k\ge 1$, the group 
$\GG_{k,a}$ is the image $\Pi_{k,a}(J_1(\ell)[\ell^k])$ of the
$\ell^k$-torsion by $\Pi_{k,a}$. It is killed by $\mgot_a^{2kg(X_1(\ell))} $ and the
restriction of $\Pi_{k,a}$ to $\GG_{k,a}$ is a bijection. 
So $\GG_a=J_1(\ell)[\mgot_a^\infty]\subset \GG$ is a divisible group. 
Let  $\eta$ be  an integer that kills $a$ in $\FF_\ell^*$ (e.g. $\eta = \ell-1$). 
Then  $\GG_{k,a}$
splits over $\FF_{p^{\eta\ell^{k+1}}}$.
\end{lemma}

The lemma below is the counterpart to 
lemma \ref{lemma:projnonsplit} in the split non-scalar case.

\begin{lemma}[The split non-scalar case]\label{lemma:projsplit}
Assume $X^2-\tau(p)X+p^{11}\pmod \ell$ has two distinct roots
$a\pmod \ell$ and $b\pmod \ell$. 
Let $k\ge 1$ be an integer and $q=p^d$ a power of $p$. Given a non
zero element in $\GG_{k,a}(\FF_{q})$, one can compute a generator  of
$V_{\ell,a}$ in  polynomial time in $\log q$, $\ell$ and $k$.
\end{lemma}

So if  $\tau(p)^2-4p^{11}$ is a non-zero  square modulo
$\ell$ we also have a simple method to construct a basis for 
the Ramanujan module $V_\ell$ modulo $p$:

We let  $a\pmod \ell$ and $b\pmod \ell$  be the two roots of
$X^2-\tau(p)X+p^{11}\pmod \ell$.
Take $q=p^{(\ell-1)\ell^4}$. We have
$\GG_a(\Fq)\supset \GG_{3,a}=\GG_{3,a}(\Fq)$ we set 
$N_q=\#J_1(\ell)(\Fq)=M_qL_q$ with $M_q$ prime to $\ell$. Let
$L_q=\ell^w$  and  
$\psi = \Pi_{w,a}\circ [M_q]$. 
The image of $J_1(\ell)(\Fq)$ by $\psi$ contains $\GG_{3,a}(\Fq)$ and is
in fact equal to  $\GG_a(\Fq)$. 
We check $\#\GG_a(\Fq)\ge \#\GG_{3,a} \ge
\ell^3$. So at least one  of the elements in $J_1(\ell)(\Fq)$ 
given by  lemma \ref{lemma:generators}
has a non-zero  image by $\psi$ for $\ell$  large enough.
We apply lemma \ref{lemma:projsplit} to this element and find a 
generator of $V_{\ell,a}$. A similar calculation produces a 
generator of  $V_{\ell,b}$. These two eigenvectors form
a basis of $V_\ell$ modulo $p$.

\medskip

All this is enough to compute the Ramanujan ideal when the Frobenius
action on it is semisimple non-scalar i.e when $\ell$ is prime to
$\tau(p)^2-4p^{11}$. 

\begin{remark}
The main simplification in this variant is that we do not need to compute pairings. 
In practice, one would just take a random degree zero $\Fq$-divisor on 
$X_1(\ell)$, multiply it by the prime to $\ell$ part of
$\#J_1(\ell)(\Fq)$   and apply  a few $B_{n,k}(T_n)$ to it.
This should usually suffice. 
\end{remark}

\begin{remark}
If $\ell$ divides $\tau(p)^2-4p^{11}$, the method described in this
section is no longer  sufficient but one can easily show that
it provides at
least one non-zero element in $V_\ell$ modulo $p$. 
\end{remark}

\section{Computing the Ramanujan subspace over $\QQ$}\label{section:lift}

Once one has
 computed the Ramanujan space $V_\ell$ inside $J_1(\ell)$ (or rather  $W_\ell$
inside  $\cJ_\ell$ the jacobian of $\cX_\ell$) modulo $p$ for many small primes $p$, one can try
to compute
this space over the rationals. 
This calculation is
described in detail in 
\cite[Section 13]{arxiv}. 
In this section we sketch a variant
of the  method presented in \cite[Section 13]{arxiv}.  We then explain
how this method should be modified to fit with the simplified
method presented in section \ref{section:semi}.
This leads us to a
 sort of generalization of the Chinese Remainder Theorem that is more adapted to
the context of polynomials with integer coefficients.

The complexity analysis of the methods presented in this section rely
on results in Arakelov theory that have been proven by 
Bas Edixhoven and  Robin de Jong,  using results by Merkl in \cite{arxiv} or J. Jorgenson and  J. Kramer 
in \cite{JJ}. In fact, the complexity analysis of the
variant  described  here requires a bit more than what has been already
given  
in \cite{arxiv}. The necessary bounds to the proof
of this variant will appear in Peter Bruin's PhD thesis \cite{bruin}.

We use the model over
$\ZZ[\frac{1}{30\ell}]$ for $\cX_\ell=X_1(5\ell)$ that is described in
section \ref{section:modularcurves2}.
 We start by fixing  a $\QQ$-rational cusp $O$ on $\cX_\ell$. This 
will be the origin of
the Jacobi map.

Let $x$ be a point in $\cJ_\ell (\bQ)$. We denote by $\theta(x)$
 the  smallest integer $k$ such that there exists
an effective divisor $D$ of degree $k$ such that  $D-kO$ belongs to  the class
 represented
by $x$ in the Picard group.
We call $\theta (x)$  the {\it stability} of $x$.
For all but finitely many primes $p$ and for any
place  $\pgot$ of $\QQ(x)$ above $p$,
one can define $\theta_\pgot (x)$ the stability
of $x$ modulo $\pgot$: the  smallest integer $k$ such that there exists
an effective divisor $D$ of degree $k$ such that  $D-kO$ belongs to  the class
 represented
by $x \bmod \pgot$ in the Picard group of $\cX_\ell\bmod p$.
We define $\theta_p(x)$ to be the minimum of all $\theta_\pgot(x)$
for all places $\pgot$ above $p$.
We note that $\theta_p(x)\le \theta_{\pgot}(x)\le \theta(x)$ whenever 
$\theta_p(x)$  is defined.
Clearly $\theta_p(x)$ is defined and equal to $\theta (x)$ for all 
large enough primes.

A  consequence of the results by Bas Edixhoven and Robin de Jong, 
extended by Peter Bruin in his forthcoming PhD thesis, see
 \cite{arxiv,bruin},  is
that, for at least half the  primes
smaller  than  $\ell^\cO$, the following holds:
 $\theta_p(x)$ is defined and equal 
to $\theta (x)$ for all  $x$  in $W_\ell$.
Notice that $\theta (x)=\theta(y)$  if $x$ and $y$ are Galois conjugate.

Now let $x$ be a non-zero point in  $W_\ell$. We can compute $x$ modulo
places $\pgot$ above 
$p$, 
 for many small (e.g. polynomial in $\ell$) primes $p$
such that $\theta_p(x)=\theta(x)$.
We only use primes such that
$\theta_p(x)=\theta(x)$ for every $x$ in $W_\ell$.

There is a unique effective divisor  $D=P_1+\dots+P_{\theta(x)}$ such that
$D-\theta(x)O$ is mapped onto $x$ by the Jacobi map. This divisor remains
unique modulo all the places  $\pgot$ in question. Further, no $P_i$
specializes to $O$ modulo any such $\pgot$.
So we choose a function $f$ on $\cX_\ell$ having no pole except at $O$.
We define e.g. 
$F(x)=f(P_1)+\cdots+f(P_{\theta(x)})$. 

We form the polynomial 

$$P_k(X)=\prod_{y\in W_\ell \mbox{ with } \theta(y)=k}(X-F(y)).$$

This polynomial has coefficients in $\QQ$.
For the above primes $p$ we have 

$$P_k(X)\bmod p=\prod_{y\in W_\ell \bmod p \mbox{ with }
  \theta_p(y)=k}(X-F(y)).$$

We set $P(X)=\prod_{k>0}P_k(X)$. If the Galois action on 
$W_\ell-\{0\}$ is transitive then $P(X)$ is likely to be irreducible and
equal to  the unique non-trivial $P_k(X)$. To be quite rigorous one should say
some  more about the choice of $f$.
See \cite[Section 22]{arxiv}.

If a reasonable  $f$  (e.g. the divisor of $f$ is $n(O-O')$ where $O'$
is another rational cusp and $n$ is the order of $O-O'$ in the jacobian) 
is
 chosen
 then Peter Bruin, improving  on  
Edixhoven, de Jong, and Merkl,  proves 
in  \cite{bruin}  that the logarithmic height of $P(X)$ is
bounded by
a polynomial in $\ell$.

If we know $W_\ell$ modulo $p$ then we can compute $P(X)$ modulo $p$ and,
provided we have  taken enough such primes $p$, we deduce $P(X)$ using
Chinese remainder theorem and the bounds proved by Edixhoven,  de Jong,
Merkl  and
Bruin.

However, if we use the simplified algorithm presented in section
\ref{section:semi} we shall only obtain $P(X)$ modulo $p$ for those
$p$
such that $\ell$ does not divide $\tau(p)^2-4p^{11}$. If $\ell$
 divides
$\tau(p)^2-4p^{11}$
then we may only obtain a non-trivial factor of $P(X)\bmod p$. This
factor  has degree
$\ell-1$ in fact.

This leads us to the following problem:

Let $P(X)$ be a degree $d\ge 2$ {\it  irreducible}\footnote{irreducible means
  here irreducible in the ring $\ZZ[X]$.} polynomial
with integer coefficients.

Let $H$ be an upper bound for the {\it naive height}
of $P(X)$: any coefficient of $P$ lies in $[-H,H]$.

Let $I$ be a positive integer and for every integer $i$ from $1$ to
$I$ assume we are given an integer  $N_i\ge 2$ and a degree $a_i$ 
{\it monic} polynomial $A_i(X)$ in $\ZZ [X]$ where $1\le a_i\le d$.
Assume the $N_i$ are pairwise coprime.

Question: assuming  $P(X)\bmod N_i$ is a multiple of 
$A_i(X)\bmod N_i$ for every  $i$, can we recover   $P(X)$,  and is
$P(X)$ the unique polynomial fulfilling all these conditions~?

We start with the following.
\begin{lemma}[Resultant and intersections]\label{lemma:resuinter}
Let $P$ and  $Q$  be two non-constant  polynomials with integer
coefficients and trivial gcd\footnote{the gcd here is the gcd in the ring
  $\ZZ[X]$.}. Let  $N\ge 2$ be an integer.
If   $P\bmod N$  and $Q\bmod N$ 
are both multiples of the same degree $d\ge 1$ monic polynomial
$A\bmod N$, then the resultant of  $P$ and  $Q$ is 
divisible by  $N^d$.
\end{lemma}
This easily follows from the resultant being given as a determinant.
\hfill $\Box$

Let  $\cP_d$ be the additive group of integer coefficient polynomials
with degree $\le d$. Let  $\rho_i : \cP_d\rightarrow
\ZZ[X]/(A_i,N_i)$ be the reduction map modulo the ideal  $(A_i,
N_i)$. 

The product map 

$$\rho = \prod_{1\le i\le I}\rho_i :
\cP_d \rightarrow \prod_{1\le i\le I}\ZZ[X]/(A_i,N_i)$$
\noindent   is surjective (Chinese remainder). 
Its kernel is therefore a lattice  $\cR$ with index  $\Theta = 
\prod_{1\le i\le I}N_i^{a_i}$ in  $\cP_d=\ZZ^{d+1}$.

If  $P_1$ and  $P_2$ are two coprime non-constant polynomials with
degree  $\le d$ and respective naive heights $K_1$ and  $K_2$,
then their resultant is bounded above by  $(2d)!K_1^{d}K_2^d$. If
further  $P_1, P_2
\in \cR$ then, according
to lemma \ref{lemma:resuinter},
$\Theta=\prod_{1\le i\le I}N_i^{a_i}$ divides the
resultant of  $P_1$ and  $P_2$. 

\begin{lemma}[Heights and intersections]
Let  $(N_i)_{1\le i\le I}$ be pairwise coprime integers.
Let $P$ be an irreducible  polynomial with integer coefficients
and degree  $d\ge 2$ and
naive height  bounded by $H$. Let  $Q$ be a polynomial with integer coefficients and
degree  $\le d$ and naive height 
bounded by $K$. Assume that for every  $i$ from  $1$ to $N$ the polynomials  $P\bmod
N_i$ and $Q\bmod N_i$ are multiples of the same monic polynomial 
 $A_i(X)\bmod N_i$  with degree  $a_i$ where  $1\le a_i\le d$. Assume
further that 

$$\prod_{1\le i\le I}N_i^{a_i} >  (2d)!H^{d}K^d.$$

Then  $Q$ is a multiple of $P$.
\end{lemma}

We observe that the  $L^2$  norm of  $P$ is  $\le H\sqrt{d+1}$. Also, if
$Q$ has $L^2$  norm  $\le H\sqrt{d+1}$ then its coefficients are 
$\le H\sqrt{d+1}$. 
Therefore if 

$$\Theta  =\prod_{1\le i\le I}N_i^{a_i} >  (2d)!(d+1)^{\frac{d}{2}}H^{2d}$$
\noindent  the polynomial  $P$ is the shortest vector in the lattice
$\cR$ for the $L^2$ norm.

Applying the  LLL algorithm to the lattice   $\cR$ we find 
(\cite[Theorem 2.6.2]{cohen})
a vector in it with $L^2$ norm $\le
2^{\frac{d}{4}}\Theta^{\frac{1}{d+1}}$. 
Taking this latter value for  $K$ we see that  if

$$ \prod_iN_i^{a_i} > (2d)!^{d+1}H^{d(d+1)}2^{\frac{d^2(d+1)}{4}}$$
\noindent then the vector output by the  LLL  algorithm is a multiple of
 $P$.

\begin{lemma}[Interpolation and lattices]
Let $d\ge 2$ be an integer. Let  $I$ be a positive integer and for every $i$ from  $1$
to  $I$ let   $N_i\ge 2$ be an integer and   $A_i(X)$  a {\it monic} polynomial
with integer coefficients and degree  $a_i$ where  $1\le
a_i\le d$. We assume the coefficients in  $A_i(X)$ lie  in the interval
 $[0,N_i[$.

We assume there exists an {\it irreducible} polynomial $P(X)$
with degree  $d$ and integer coefficients and naive height  $\le H$ such
that   $P(X)\bmod N_i$ is a multiple of 
$A_i(X)\bmod N_i$ for all $i$.

We assume the  $N_i$ are pairwise coprime and 
$$\prod_{1\le i\le I}N_i^{a_i}  > (2d)!^{d+1}H^{d(d+1)}2^{\frac{d^2(d+1)}{4}}.$$

Then  $P(X)$ is  the unique polynomial fulfilling all these conditions
and it can be computed  from the $(N_i,A_i(X))$
by a  deterministic Turing machine in time polynomial  in  $d$,  $\log
H$ and  $I$,  and the $\log N_i$. 
\end{lemma}

Note that the dependency on  $I$ and  $\log N_i$  is harmless 
because one may remove some information if there is too much of it.
We can always do with some   $I$ and  $\log N_i$  that are polynomial
in $d$ and  $\log H$.

This lemma shows that we can compute (lift) the Ramanujan module
$W_\ell$ using  the
simplified algorithm of section \ref{section:semi},
even if the action of the Frobenius  at $p$ on $W_\ell$ is not
semisimple for any auxiliary prime  $p$. 

\section{Are there many semi simple  pairs $(\ell,p)$ ?}\label{section:elem}

We have seen in section \ref{section:semi} that the computation of
$V_\ell$ modulo $p$  becomes
simpler
whenever the two primes $p$ and $\ell$ satisfy the condition that 
$\ell$  is prime to $\tau(p)^2-4p^{11}$. If this is the case, we say
that
the pair $(\ell,p)$ is good (otherwise it is bad).

In the situation of section \ref{section:lift} we are given a fixed
prime $\ell$ and we look for primes $p$ such that $(\ell,p)$ is
good. We need these primes $p$ to be bounded by a polynomial in
$\ell$. And there should be enough of them that we can find them by 
random search.

This leads us to the following definition.

\begin{definition}[What bad and good means in this section]\label{definition:goodprimes}
We say that a pair $(\ell,p)$ of prime integers is bad  if  $\ell$ divides
$\tau(p)^2-4p^{11}$. Otherwise it is good. Let  $c>1$ be a real.   We say that a given  prime
$\ell$ is $c$-bad if $(\ell,p)$ is bad for at  least half the 
primes $p\le \ell^c$. Otherwise it is $c$-good.
\end{definition}

In this section we give an elementary unconditional proof that there
are enough good primes $\ell$.
Let $\alpha$, $\beta$, $\gamma$ and $\delta$ be four positive
constants such that for every   integer $k\ge 2$ the $k$-th prime
$p_k$ satisfies $\alpha k\log k \le p_k \le \beta k\log k$ and for
every real 
$x\ge 2$ the
arithmetic function $\pi(x)$ giving the number of primes $\le x$  satisfies $\gamma x(\log x)^{-1} \le \pi (x)
\le \delta x (\log x)^{-1}$.

Work by Tchebitchef  allows $\gamma =\frac{1}{3}$ and $\delta =
\frac{5}{4}$. Work by Rosser \cite{rosser} shows that $\alpha = 1$ is fine. Rosser
also proved that $p_k\le k(\log k+\log\log k)$ for $k\ge 6$. So we can
take $\beta=2.17$ for example. I thank Guillaume Hanrot for pointing out 
these references to me.

Let $X\ge 3$ be an integer. Let   $L$  be 
the $X$-th prime integer.  Let $\cX(c,X)$ be the set of pairs of primes $(\ell,p)$ with $\ell
\le L$ and $p\le \ell^c$.
We set $\ell_1=p_1=2$, $\ell_2=p_2=3$, \ldots 
 the successive prime integers.
 Let $P$ be the largest prime 
$\le L^c$ and let  $Y$ be the integer such that $P=p_Y$.
One has $L\le \beta X\log X$ and $P\le \beta^cX^c(\log X)^c$ and $Y\le
P$.

Since $\tau(p)^2-4p^{11}$ has at most $\log_2(4p^{11})$ prime divisors, 
there are at most $Y(2+11\log_2 P)$ bad pairs and this is $\le
51c\beta^cX^c(\log X)^{c+1}$ provided $X\ge \beta$.  We want to bound from above  the number of bad $\ell \le L$. 
The worst case is when the smallest $\ell$ are bad. Assume all primes
$\ell \le \ell_x$ are bad. The number of bad pairs is then  at
least 

$$\frac{1}{2}\sum_{1\le
  k\le x} \pi(\ell_k^c)\ge \frac{\gamma\alpha^c}{2}\sum_{\frac{3}{\alpha} \le k\le x}
\frac{k^c(\log k)^c}{c\log \alpha + c\log k + c\log\log k}\ge
\frac{\gamma\alpha^c}{4c} \sum_{\frac{3}{\alpha} \le k\le x} k^c(\log
k )^{c-1}$$
\noindent and this is at least

$$\frac{\gamma\alpha^c}{4c(c+1)}{\left(x^{c+1} -\left(
  \frac{3}{\alpha}\right)^{c+1}\right)}\ge \frac{\gamma
  \alpha^c}{8c(c+1)}x^{c+1}$$
\noindent   provided $x\ge 6/\alpha$. Assume
at  least half of the primes $\ell\le L$ are bad. Then the number of bad
pairs
is at least $\frac{\gamma   \alpha^c}{8c(c+1)}(X/2)^{c+1}$  provided
$X\ge 12/\alpha$. So 

$$\frac{\gamma   \alpha^c}{8c(c+1)}(X/2)^{c+1}\le 51c\beta^cX^c(\log
X)^{c+1}$$
\noindent so 
$$\frac{X}{(\log X)^{c+1}}\le 816\left(\frac{2\beta}{\alpha}
\right)^cc^2(c+1)\gamma^{-1}.$$

We call  $a$ the right-hand side in the above inequality. We 
set $Z=X^{\frac{1}{c+1}}$
and we have
$\frac{Z}{\log Z}\le (c+1)a^{\frac{1}{c+1}}$.
Since  $\log Z\le \sqrt Z$ we have $Z\le
(c+1)^2a^{\frac{2}{c+1}}$ and $X\le (c+1)^{2(c+1)}a^{2}$.

\begin{lemma}\label{lemma:goodprimes}
Let $\alpha$, $\beta$, $\gamma$ and $\delta$ be the four constants
 introduced before definition \ref{definition:goodprimes} 
 above. Let $c>1$ be a real number. Assume $X$ is an integer  bigger
 than   $816^2c^4(c+1)^{2(c+2)}\left(\frac{2\beta}{\alpha}
\right)^{2c}\gamma^{-2}$. Then at least half among  the $X$ first
primes are $c$-good.
\end{lemma}

\begin{lemma}[Effective  bound for the density of good primes $\ell$]\label{lemma:goodprimes2}
Let $c>1$ be a real number. Assume $X$ is an integer  bigger
 than $2^{23+5c}c^4(c+1)^{2(c+2)}$. Then at least half among  the $X$ first
primes are $c$-good.
\end{lemma}

\appendix

\section{A GP-PARI code for  Puiseux
expansions at singular branches of modular curves}\label{app:A}

Below are a few lines of GP-PARI code (see \cite{pari}) that compute
the  expansions of $x_{\alpha,\beta}$ as series in $b^{-\frac{1}{\ell}}$
with coefficients in a finite field containing a primitive 
$\ell$-th root of unity. We use the methods and notation
given
in  section
\ref{section:modularcurves2}, before the statement of lemma \ref{lemma:computingCl2}.

Our  code 
computes the $q$-series for the modular function $j$ as 

$$j(q)=1728E_4^3(q)(E_4^3(q)-E_6^2(q))^{-1}$$
\noindent  where

$$E_4(q)=1+240\sum_{n\ge 1} \frac{n^3q^n}{1-q^n}$$
\noindent
and $$E_6(q)=1-504\sum_{n\ge 1} \frac{n^5q^n}{1-q^n}.$$

The expansions for the $x_{\alpha,\beta}$ are then obtained through
standard operations on series like product, sum, reversion,
composition.

\begin{verbatim}
{ser(aa,bb,prec,ell,p,z,b,jc,E4,E6,D,jq,qc,gc,w,x)=
ell=7;
p=953;
z=Mod(431,p);
b=1/c;
jc=(b^4-12*b^3+14*b^2+12*b+1)^3/b^5/(b^2-11*b-1);
E4=sum(n=1,prec, n^3*q^n/(1-q^n))*240+1+O(q^prec);
E6=sum(n=1,prec, -n^5*q^n/(1-q^n))*504+1+O(q^prec);
D=(E4^3-E6^2)/1728;
jq=E4^3/D;
qc=subst(serreverse(1/jq),q,1/jc+O(c^prec));
gc= -36*b*(b^2-11*b-1)*deriv(qc)*(-c^2)/5/qc;
w=z^aa*Q^(2+5*bb);
xabs=Mod(1,p)*(1/12
+sum(n=1,prec,
w*Q^(5*ell*n)/(1-w*Q^(5*ell*n))^2+O(Q^(5*ell*prec)))
+w/(1-w)^2
+sum(n=1,prec,
Q^(5*ell*n)/w/(1-(w)^(-1)*Q^(5*ell*n))^2+O(Q^(5*ell*prec)))
-2*sum(n=1,prec,
n*Q^(5*ell*n)/(1-Q^(5*ell*n))+O(Q^(5*ell*prec )) ));
cQ=subst(serreverse((qc/c^5)^(1/5)*c),c,Q^ell);
bQ=1/cQ;
gQ=subst(gc,c,cQ);
XabQ=(gQ*xabs-3*(bQ^2-6*bQ+1) )/36;
QC=subst(serreverse(1/((bQ*Q^ell)^(1/ell)/Q)),Q,C);
XabC=subst(XabQ,Q,QC);
}
\end{verbatim}


\section{A Magma code that computes the zeta function
of  modular curves}\label{app:B}

Below are a few lines  written in the Magma language (see \cite{magma}). They compute
the  characteristic polynomial
of the Frobenius of $X_1(5\ell)/\Fp$ using the methods given in the proof of
lemma
\ref{lemma:manin2}. 

\begin{verbatim}
ZZ:=IntegerRing();
l:=11;
N:=5*11;
QN:=CyclotomicField(EulerPhi(N));
R1<T>:=PolynomialRing(QN,1);
R2<T,U>:=PolynomialRing(QN,2);
G := DirichletGroup(N,QN);
chars := Elements(G);
gen4:=chars[2];
gen10:=chars[5];
Genus(Gamma1(N));
charsmc:=[gen4,gen4^2,gen4^4, gen4*gen10,gen4^2*gen10,
gen10,gen4*gen10^2,gen4^2*gen10^2,gen10^2 , gen4*gen10^5,
gen4^2*gen10^5,gen10^5];
p:=101;
PT:= R2 ! 1;
W:=1;
g:=1;

for eps in charsmc do

M := ModularForms([eps],2); 
P:= R2 ! Evaluate(HeckePolynomial(CuspidalSubspace(M),p),T);
g:=Degree(P,T);
W :=   Evaluate(P,[ T+Evaluate(eps,p)*p/T,  1])*T^g;
PT:=PT*W;

end for;

PT := R2 ! PT;

k:=2;
PTk:= Resultant(PT, T^k-U,T);
\end{verbatim}

\end{document}